\documentclass[11pt]{article}

\usepackage{amstext,amssymb,amsmath,amsbsy}
\usepackage{hyperref}
\usepackage{amscd}
\usepackage{amsfonts}
\usepackage{indentfirst}
\usepackage{verbatim}
\usepackage{amsmath}
\usepackage{amsthm}
\usepackage{enumerate}
\usepackage{graphicx}
\usepackage{color, soul}
\usepackage[OT1]{fontenc}
\usepackage[latin1]{inputenc}
\usepackage[english]{babel}
\usepackage{amssymb}
\usepackage{subfig}
\usepackage{algorithm}
\usepackage{algpseudocode}

\usepackage{mathtools}

\newcommand{\R}{\mathbb{R}}
\newcommand{\N}{\mathbb{N}}

\newcommand{\x}{{\bf x}}

\newcommand{\p}{{\bf p}}

\newcommand{\Div}{{\rm div}}

\newtheorem{Theorem}{Theorem}[section]
\newtheorem{Lemma}{Lemma}[section]

\newtheorem{Corollary}{Corollary}[section]
\newtheorem{remark}{Remark}[section]

\newtheorem*{Assumption*}{Assumption}

\newtheorem{problem}{Problem}[section]
\newtheorem*{problem*}{Problem}
\numberwithin{equation}{section}

\begin{document}

\title{The Carleman contraction mapping method for quasilinear elliptic equations with over-determined boundary data}

\author{Loc H. Nguyen\thanks{Department of Mathematics and Statistics, University of North Carolina at
Charlotte, Charlotte, NC, 28223, USA, \texttt{loc.nguyen@uncc.edu}.}}  


\date{}
\maketitle
\begin{abstract}
	We propose a globally convergent numerical method to compute solutions to a general class of quasi-linear PDEs with both Neumann and Dirichlet boundary conditions. 
	Combining the quasi-reversibility method and a suitable Carleman weight function, we define a map of which fixed point is the solution to the PDE under consideration.
	To find this fixed point, we define a recursive sequence with an arbitrary initial term using the same manner as in the proof of the contraction principle.
	Applying a Carleman estimate, we show that the sequence above converges to the desired solution.
	On the other hand, we also show that our method delivers reliable solutions even when the given data are noisy.
	Numerical examples are presented.
\end{abstract}




\section{Introduction}\label{sec1}

Let $\Omega$ be an open and bounded domain in $\mathbb{R}^d$, $d \geq 2$, with a smooth boundary. 
Let $F: \overline \Omega \times \mathbb{R} \times \mathbb{R}^d \to \mathbb{R}$ be a real-valued function in the class $C^2$.
Let $A$ be a matrix-valued function $\overline \Omega \to \mathbb{R}^{d \times d}$
satisfying
\begin{enumerate}
    \item $A$ is in the class $C^2(\overline \Omega, \mathbb{R}^{d \times d} )$,
    \item $A$ is symmetric; i.e. $A^{\rm T} = A$,
    \item there are  positive constants $\Lambda_1$ and $\Lambda_2$ such that
    \begin{equation*}
        \Lambda_1 \vert \xi\vert^2 \leq A(\x)\xi \cdot \xi \leq \Lambda_2 \vert \xi \vert^2
        \quad \mbox{for all } \x \in \overline \Omega, \xi \in \mathbb{R}^d.
    \end{equation*}
\end{enumerate}
Consider the over-determined boundary value problem
\begin{equation}
    \left\{
        \begin{array}{ll}
             {\rm Div} (A(\x) \nabla u(\x))+ F(\x, u(\x), \nabla u(\x)) = 0 &\x \in \Omega,  \\
             u(\x) = f(\x) &\x \in \partial \Omega,\\
             A(\x)\nabla u(\x) \cdot \nu(\x) = g(\x) &\x \in \partial \Omega
        \end{array}
    \right.
    \label{main_eqn}
\end{equation}
where $f$ and $g$ are two given functions, which are the noisy measurements in some applied contexts; e.g., see the inverse problem in Section \ref{inverse1}.
The main aim of this paper is to develop a numerical method to solve
the following problem.
\begin{problem}
    Let $f^*$ and $g^*$ be the noiseless versions of $f$ and $g$ respectively. 
    Assume that problem  
    \begin{equation}
    \left\{
        \begin{array}{ll}
             {\rm Div} (A(\x) \nabla u^*(\x))+ F(\x, u^*(\x), \nabla u^*(\x)) = 0 &\x \in \Omega,  \\
             u^*(\x) = f^*(\x) &\x \in \partial \Omega,\\
             A(\x)\nabla u^*(\x) \cdot \nu(\x) = g^*(\x) &\x \in \partial \Omega
        \end{array}
    \right.
    \label{main_eqn_exact}
\end{equation}
has a unique solution $u^*$. Given the noisy data $f$ and $g$, compute an approximation of $u^*$.
\label{p}
\end{problem}

In the statement of Problem \ref{p}, we request both Dirichlet and Neumann boundary data. In the theory of PDEs, one of these data might be sufficient to determine the solution uniquely.
Hence, Problem \ref{p} is over-determined.
Our study accepts this redundant weakness because we solve Problem \ref{p} for the needs of inverse problems.
The application in inverse problems is explained as follows. 
Recently, we numerically solved  several inverse problems by a unified framework, see e.g., \cite{VoKlibanovNguyen:IP2020, Khoaelal:IPSE2021, KhoaKlibanovLoc:SIAMImaging2020, KlibanovLeNguyenIPI2021, KlibanovNguyen:ip2019, LeNguyen:jiip2022, Nguyen:CAMWA2020, Nguyens:jiip2020, KlibanovAlexeyNguyen:SISC2019}.
This framework has two steps.
\begin{itemize}
\item In step 1, we introduce a change of variable to reduce the given inverse problem to a system of quasi-linear PDEs with Cauchy boundary data.
\item In step 2, we solve the over-determined system obtained in Step 1. The computed solution yields the solution to the inverse problem under consideration.
\end{itemize}  
The goal of Problem \ref{p} is to address step 2 above. 
That is how to solve a system of quasi-linear elliptic equations with Dirichlet and Neumann data.
For convincing purposes, we solve an inverse source problem for a nonlinear model in Section \ref{inverse1}. 
This serves as an example of reducing a challenging nonlinear inverse problem to a system of PDEs of which this kind of over-determined data is available.
We also cite to \cite{KlibanovLeNguyenIPI2021, LeKlibanov:ip2022} for using this framework to solve the inverse scattering problem in the time domain with experimental data.
For simplicity, we solve a single equation rather than solving a system of quasi-linear elliptic PDEs with Cauchy boundary data. This simplification does not weaken the paper because our analysis and numerical implementation can be directly extended for systems of quasi-linear equations.

As mentioned in the paragraph above, in the theory of PDEs, one might need only one boundary condition to determine the solution to \eqref{main_eqn}.
However, this might not be true in some specific circumstances. 
For example, the equation $y''(t) + \pi^2y(t) = 0$, $t \in (0, 1)$, with $y(0) = y(1) = 0$ has infinitely many solutions  $y(t) = C \sin(\pi  t)$, $C \in \R$. 
We cite to \cite{DancerSchmitt:ops1987, DucLocLuc:ejde2013, LocSchmitt:die2009} for more examples in which quasi-linear elliptic PDEs, with one boundary condition, have multiple solutions.
Therefore, our drawback when requiring over-determined data might be acceptable.
On the other hand, up to the author's knowledge, a numerical method to solve quasi-linear elliptic PDE with only one boundary condition is not yet developed unless more information about the solutions is known.

Since \eqref{main_eqn} involves both Dirichlet and Neumann conditions, \eqref{main_eqn} might not have a solution; especially, when the measured data $f$ and $g$ contain significant noise. 
Computing the solution to \eqref{main_eqn} might be impossible.
In this case, we understand the solution to \eqref{main_eqn} as the limit of a sequence obtained by iteratively solving linear least squares optimization problems.
  Assuming that \eqref{main_eqn} with noiseless data has a unique smooth solution $u^*$ and given noisy data, we will rigorously prove that this sequence approximates $u^*$.
  This result is one of the most critical points of the current paper.

Our proposed numerical method to solve quasi-linear elliptic equations with Cauchy data in this paper has two crucial features: fast and global. 
By ``fast", we mean that the method converges at the exponential rate with respect to the number of iterations.
By ``global", we mean that our method does not require a good initial guess of the true solution to the problem under consideration.
Both features are the crucial strengths of this paper since it is well-known that the widely used optimization-based methods for solving nonlinear equations are local and time-consuming.

 In the statement of Problem \ref{p}, we have imposed a condition about the existence and uniqueness of bounded solutions to \eqref{main_eqn}, with $f$ and $g$ replaced by $f^*$ and $g^*$ respectively. 
This condition can be interpreted as follows. Our target is to provide a new tool to solve nonlinear inverse problems using the framework mentioned above. 
In these applications, the solution $u$ represents some physical quantities related to heat distribution or wave propagation; see Section \ref{inverse1} for an example. So, when the measured data $f$ and $g$ are perfectly noiseless, \eqref{main_eqn} has a solution that is such a physical quantity. 
So, the existence is clear from the physical point of view. The uniqueness is due to the presence of both Dirichlet and Neumann conditions and the unique continuation principle.

Problem \ref{p} is exciting and challenging partly because our target is to compute $u^*$ when the noisy data $f$ and $g$ are given while the corresponding noiseless ones $f^*$ and $g^*$ are unknown.
A natural approach to compute the solution to \eqref{main_eqn} is to minimize a least squares functional.
A typical example of such a functional is
\begin{multline}
     H^p(\Omega) \ni u \mapsto 
    J(u) =  \int_{\Omega}\big\vert{\rm Div} (A(\x) \nabla u(\x))+ F(\x, u(\x), \nabla u(\x))\big\vert^2d\x
     \\
     + \mbox{a regularization term}
     \label{3}
\end{multline}
subject to the boundary conditions in \eqref{main_eqn}.
One takes the minimizer of the functional $J$ in \eqref{3} as a solution to Problem \ref{p}. This approach is based on optimization. 
It has three main drawbacks:
\begin{enumerate}
    \item \label{d1} $J$ might be nonconvex. It might have multiple local minimizers. Therefore, a good initial guess of the true solution $u^*$ is required.
    \item \label{d2} The computational cost is  expensive.
    \item \label{d3} It is not clear that the minimizer is an approximation of $u^*.$
\end{enumerate}
Recently, we have developed the convexification method, see \cite{LeNguyen:JSC2022}, and the Carleman weighted linearization method, see \cite{LeNguyenTran:CAMWA2022}, to solve Problem \ref{p}.
\begin{itemize}
    \item 
The key point of the convexification method is to include suitable Carleman weight functions into the formulation of the mismatch functional $J$. 
By using  Carleman estimates, one can prove that the new mismatch functional is strictly convex. One also can prove that the unique minimizer is a good approximation of $u^*$. 
Hence, drawbacks \ref{d1} and \ref{d3} can be overcome.
The convexification method was first introduced in \cite{KlibanovIoussoupova:SMA1995} and then was developed intensively.
We refer the reader to \cite{Klibanov:ip2015,Klibanov:sjma1997,Klibanov:nw1997,KlibanovNik:ra2017,KlibanovKolesov:cma2019,KlibanovLiZhang:ip2019,KhoaKlibanovLoc:SIAMImaging2020,KlibanovLiZhang:SIAM2019, KlibanovLeNguyenIPI2021, LeNguyen:JSC2022, KlibanovLiBook}
for some important works in this area and their applications to solving a
variety kinds of inverse problems.
However, the computation due to the convexification method is time-consuming.
\item 
We have introduced in \cite{LeNguyenTran:CAMWA2022} another method, also based on Carleman estimates, to solve Problem \ref{p}. 
The method in \cite{LeNguyenTran:CAMWA2022}  is inspired by Carleman estimates and linearization similar to the Newton method. We have shown in \cite{LeNguyenTran:CAMWA2022} that the combination of Carleman estimates and linearization allows us to compute $u^*$ quickly without requesting a good initial guess.
The Carleman-Newton method successfully solved a nonlinear inverse problem in \cite{AbhishekLeNguyenKhan} and computed numerical solutions to Hamilton-Jacobi equations in \cite{LeNguyenTran:CAMWA2022}.
\end{itemize}

The contribution of this paper is to introduce another globally convergent numerical method based on a Carleman estimate and the classical contraction principle.
More precisely, our approach is first to define a map $\Phi$ such that the desired solution is the fixed point of this map. The construction of $\Phi$ combines the Carleman weight function and the quasi-reversibility method to solve over-determined linear PDEs (see \cite{LattesLions:e1969} for the original work for the quasi-reversibility method).
Using a suitable Carleman estimate, we rigorously prove that  $\Phi$ is a contraction map. This leads to a numerical method to solve Problem \ref{p}. We simply approximate the desired solution by $u_n = \Phi^n(u_0)$ where $\Phi^n = \Phi \circ \Phi \dots \circ \Phi$ ($n$ times) and $u_0$ is an arbitrary function.
The main theorems in this paper confirm that our function $\Phi$ is a contraction map and that the sequence $\{u_n\}_{n \geq 0}$ above converges to the true solution. Imposing some technical assumptions, we will prove that the stability with respect to the noise contained in the given data is of the Lipschitz type.
We also refer to \cite{ BAUDOUIN:SIAMNumAna:2017, Baudouin:SIAM2021, LeNguyen:jiip2022} for similar works for the case when the data has no noise and refer to \cite{LeKlibanov:ip2022, NguyenKlibanov:ip2022} for the proof of a similar result for hyperbolic equations.
The strengths of our new approach include the fact that
\begin{enumerate}
	\item it does not require a good initial guess;
	\item it is quite general in the sense that no special structure is imposed on the nonlinearity $F$; 
	\item the convergence rate is $O(\theta^n)$ where $\theta \in (0, 1)$ and $n$ is the number of iterations.
\end{enumerate}

The paper is organized as follows.
In Section \ref{inverse1}, we present an inverse source problem that motivates the study of Problem \ref{p}.
In Section \ref{sec2}, we introduce the contraction map $\Phi$, which plays the key role in solving Problem \ref{p}.
In Section \ref{sec3}, we show that the fixed point of $\Phi$ is an approximation of the solution to Problem \ref{p}. We investigate the hehavior of this approximation as the noise in the boundary data tends to $0$. 
Section \ref{sec4} is for the numerical study.
Section \ref{sec6} is for some concluding remarks.



\section{An inverse source problem for nonlinear hyperbolic equations} \label{inverse1}

In this section, we provide an example from which
 Problem \ref{p} arises.
 Let $T > 0$ represent the final time and $G: \R^d \times \R \times \R^d \to \R$ be a smooth function.
Consider the wave function $w: \R^d \times (0, T) \to \R$ satisfying the following initial value problem
\begin{equation}
    \left\{
        \begin{array}{ll}
              w_{tt}(\x, t) = \Delta w(\x, t) + G(\x, w(\x, t), \nabla w(\x, t))&
             \x \in \R^d, t \in (0, T)\\
             w(\x, 0) = p(\x)& \x \in \R^d,\\
             w_t(\x, 0) = 0 &\x \in \R^d.
        \end{array}
    \right.
    \label{nonlinearwave}
\end{equation}
 Here, $p(\x)$ is a source term that generates the wave. 
 The nonlinear inverse source problem we are interested in is formulated as follows.
\begin{problem}[Inverse Source Problem for hyperbolic equations]
Assume that $p$ is compactly supported in a smooth and bounded domain $\Omega$. 
Compute the source function $p(\x)$, $\x \in \Omega$, from the measurements of 
\begin{equation}
	f_1(\x, t) = w(\x, t) \quad
	\mbox{and}
	\quad
	f_2(\x, t) = \partial_{\nu}w(\x, t) 
	\label{data_isp}
\end{equation}
for all $\x \in \partial \Omega$, $t \in [0, T].$
\label{ISP}
\end{problem}
Problem \ref{ISP} can be considered the nonlinear version of the thermo/photo-acoustics tomography problem arising from bio-medical imaging. 
The experiment leading to this problem is as follows, see \cite{Krugeretal:mp1995, Krugerelal:mp1999, Oraevskyelal:ps1994}.
One sends non-ionizing laser pulses or microwaves to a biological tissue under inspection (for instance, a woman's breast in mammography). Some energy will be absorbed and converted into heat, causing a thermal expansion and a subsequence ultrasonic wave propagating in space. The ultrasonic pressures on a surface around the tissue are measured.  Finding some initial information about the pressures from these measurements yields the structure inside this tissue.
Most works in the field of thermo/photo-acoustics tomography address the problem when the governing hyperbolic equation is linear, while the study for nonlinear cases is very limited \cite{NguyenKlibanov:ip2022}. 
We list here some widely used methods for the linear models.
In the case when the waves propagate in the free space, one can find explicit reconstruction formulas in \cite{DoKunyansky:ip2018,  Haltmeier:cma2013, Natterer:ipi2012, Linh:ipi2009}, the time reversal method \cite{ KatsnelsonNguyen:aml2018, Hristova:ip2009, HristovaKuchmentLinh:ip2006, Stefanov:ip2009, Stefanov:ip2011}, the quasi-reversibility method \cite{ClasonKlibanov:sjsc2007, LeNguyenNguyenPowell:JOSC2021} and the iterative methods \cite{Huangetal:IEEE2013, Paltaufetal:ip2007, Paltaufetal:osa2002}.
The publications above study thermo/photo-acoustics tomography for simple models for non-damping and isotropic media. 
The reader can find publications about thermo/photo-acoustics tomography for more complicated model involving a damping term or attenuation term \cite{Ammarielal:sp2012, Ammarietal:cm2011, Haltmeier:jmiv2019, Acosta:jde2018, Burgholzer:pspie2007, Homan:ipi2013, Kowar:SISI2014, Kowar:sp2012, Nachman1990}.
In this section, we propose another method based on our solver of Problem \ref{p}.

Let $\{\Psi_n\}_{n \geq 1}$ be the orthonormal basis of $L^2(0, T)$ originally introduced in \cite{Klibanov:jiip2017} 
and define
\begin{equation}
	w_n(\x) = \int_{0}^{T} w(\x,t) \Psi_n(t) \,dt
	\quad \text{ for } n \geq 1, \x \in \Omega.
	\label{zn}
\end{equation}
In computation,  we can approximate
\begin{equation}
w(\x,t) =\sum_{n=1}^\infty w_n(\x) \Psi_n(t) \approx \sum_{n=1}^N w_n(\x) \Psi_n(t),
\label{wt}
\end{equation}
for a suitable cut-off number $N\in \N$.
Then, due to the governing equation in \eqref{nonlinearwave},
the vector $W_N = (w_1, w_2, \dots, w_N)$ ``approximately" satisfies 
\begin{equation}
     \sum_{n=1}^N w_n(\x) \Psi_n''(t) 
     =  \sum_{n=1}^N \Delta w_n(\x) \Psi_n(t) 
     + G\big(\x, \sum_{i=1}^N w_n(\x) \Psi_n(t),  \sum_{n=1}^N \nabla w_i(\x) \Psi_n(t)\big)
     \label{10_10}
\end{equation}
for all $\x \in \Omega$ and $t \in (0, T)$.
For each $m \in \{1, \dots, N\},$ we multiply $\Psi_m(t)$ to both sides of \eqref{10_10} and then integrate the resulting equation.
We obtain
\begin{multline}
     \sum_{n=1}^N w_n(\x) \int_{0}^T\Psi_n''(t)  \Psi_m(t) dt
     =  \sum_{n=1}^N \Delta w_n(\x) \Psi_n(t) \Psi_m(t)dt
     \\
     + \int_{0}^T G\Big(\x, \sum_{n=1}^N w_n(\x) \Psi_n(t),  \sum_{n=1}^N \nabla w_n(\x) \Psi_n(t)\Big) \Psi_m(t)dt
     \label{11_11}
\end{multline}
for all $\x \in \Omega.$
Denote
\begin{align*}
    &W(\x) = (w_1(\x), w_2(\x), \dots, w_N(\x))^{\rm T},\\
    &S = (s_{mn})_{m, n = 1}^N,\\
    &\mathcal{G}(\x, W(\x), \nabla W(\x)) = (g_1(\x, W(\x), \nabla W(\x)),  \dots, g_N(\x, W(\x), \nabla W(\x)))^{\rm T}
\end{align*}
where
\begin{equation*}
    s_{mn} = \int_{0}^T\Psi_n''(t)  \Psi_m(t) dt
\end{equation*}
and
\begin{equation*}
    g_m(\x, W(\x), \nabla W(\x)) = \int_{0}^T G\Big(\x, \sum_{n=1}^N w_n(\x) \Psi_n(t),  \sum_{n=1}^N \nabla w_n(\x) \Psi_n(t)\Big) \Psi_m(t)dt
\end{equation*}
for all $\x \in \Omega.$
    It follows from \eqref{11_11} that
    \begin{equation}
        \Delta W(\x) + \mathcal F(\x, W(\x), \nabla W(\x)) = 0
        \quad \mbox{for all }
        \x \in \Omega
            \label{eqn_isp}
    \end{equation}
where
\begin{equation*}
    \mathcal F(\x, W(\x), \nabla W(\x)) = \mathcal G(\x, W(\x), \nabla W(\x)) - SW(\x).
\end{equation*}

Boundary conditions for the vector-valued function $W$ can be computed from the given boundary data in the statement of Problem \ref{ISP}.
Its follows from \eqref{data_isp} and \eqref{zn} that for all $\x \in \partial \Omega$
\begin{equation}
    W(\x) = {\bf f}(\x) =\Big(\int_0^T f_1(\x, t) \Psi_n(t)dt\Big)_{n = 1}^N
    \label{dir_isp}
\end{equation}
and 
\begin{equation}
    \partial_{\nu} W(\x) = {\bf g}(\x) = \Big( \int_0^T f_2(\x, t) \Psi_n(t)dt\Big)_{n = 1}^N.
    \label{neu_isp}
\end{equation}
\begin{remark}
Computing a function $W = W^{\rm comp} = (w_1^{\rm comp}, \dots, w_N^{\rm comp})^{\rm T}$ from \eqref{eqn_isp}, \eqref{dir_isp} and \eqref{neu_isp} is the goal of Problem \ref{p}. 
This partly shows the significance of the study in this paper.
\end{remark}
Having $W^{\rm comp}$ in hand, due to \eqref{wt}, we can evaluate the source term $p$ via
\begin{equation}
    p(\x) = \sum_{n = 1}^N w_n^{\rm comp}(\x) \Psi_n(0)
    \label{p_comp}
\end{equation}
for all $\x \in \Omega.$
Since solving inverse problems is out of the scope of this paper, we will present the details about this method for Problem \ref{ISP} and some numerical results in a near future publication.

\section{The Carleman contraction principle} \label{sec2}

In this section, we establish a Carleman contraction method to solve quasi-linear PDEs.
The main tool that guarantees the success of our method is a Carleman estimate.
Carleman estimates are great tools in the study of PDEs. They were first used to prove the unique continuation principle, see, e.g., \cite{Carleman:1933, Protter:1960AMS}.
The use of Carleman estimates quickly became a powerful tool in many areas of PDEs, especially in both theoretical and numerical methods for inverse problems, see, e.g., \cite{BukhgeimKlibanov:smd1981, BeilinaKlibanovBook, KlibanovLiBook, VoKlibanovNguyen:IP2020, KhoaKlibanovLoc:SIAMImaging2020, KlibanovNguyen:ip2019, LeNguyenNguyenPowell:JOSC2021, LocNguyen:ip2019}.
Carleman estimates were used in cloaking \cite{MinhLoc:tams2015} and in the area of computing solution to Hamilton-Jacobi equations \cite{KlibanovNguyenTran:JCP2022, LeNguyenTran:CAMWA2022}.
We recall a useful Carleman estimate which is important for us in the proof of the main theorem in this paper.
Let $\x_0$ be a point in $\mathbb{R}^d \setminus \overline \Omega$ such that $r(\x) = \vert\x  -\x_0\vert > 1$ for all $\x \in \Omega.$
For each $\beta > 0$, define 
\begin{equation}
	\mu_\beta(\x) = r^{-\beta}(\x)= \vert\x - \x_0\vert^{-\beta}  \quad \mbox{for all } \x \in \overline \Omega.
	\label{mu}
\end{equation}

We have the following lemma.

\begin{Lemma}[Carleman estimate]\label{lem:Carleman}
There exist  positive constants $\beta_0$ depending only on $\x_0$, $\Omega$, $\Lambda,$  and $d$ such that for all function $v \in C^2(\overline \Omega)$ satisfying
	\begin{equation}
		v(\x) = \partial_{\nu} v(\x) = 0 \quad \mbox{for all } 
		\x \in \partial \Omega,
		\label{3.1}
	\end{equation}
	the following estimate holds true
	\begin{equation}
		\int_{\Omega} e^{2\lambda \mu_\beta(\x)}\vert{\rm Div}(A \nabla v) \vert^2d\x
		\geq
		 C\lambda  \int_{\Omega}  e^{2\lambda \mu_\beta(\x)}\vert\nabla v(\x)\vert^2\,d\x
		+ C\lambda^3  \int_{\Omega}   e^{2\lambda \mu_\beta(\x)}\vert v(\x)\vert^2\,d\x	
		\label{Car est}	
	\end{equation}
	for all $\beta \geq \beta_0$ and $\lambda \geq \lambda_0$. 
	Here, $\lambda_0 = \lambda_0( \x_0,  \Omega, A, d, \beta)$ and $C = C( \x_0,  \Omega, A, d, \beta) > 0$ depend only on the listed parameters.
	\label{lemma carl}
\end{Lemma}

Lemma \ref{lemma carl} is a direct consequence of \cite[ Lemma 5]{MinhLoc:tams2015}. We refer the reader to \cite[Lemma 2.1]{LeNguyenTran:CAMWA2022} for details of the proof.
An alternative way to obtain \eqref{Car est}, with another Carleman weight function, is to apply the Carleman estimate in  \cite[Chapter 4, Section 1, Lemma 3]{Lavrentiev:AMS1986} for general parabolic operators. 
The arguments to obtain \eqref{Car est} using \cite[Chapter 4, Section 1, Lemma 3]{Lavrentiev:AMS1986} are similar to that in \cite[Section 3]{LeNguyenNguyenPowell:JOSC2021} with the Laplacian replaced by the operator ${\rm Div} (A\nabla \cdot)$. 
	We especially draw the reader's attention to different forms of Carleman estimates for all three kinds of differential operators (elliptic, parabolic, and hyperbolic) and their applications in inverse problems and computational mathematics \cite{BeilinaKlibanovBook, BukhgeimKlibanov:smd1981, KlibanovLiBook, LocNguyen:ip2019}. 
	It is worth mentioning that some Carleman estimates hold true for all functions $v$ satisfying $v\vert_{\partial \Omega} = 0$ and $\partial_{\nu} v\vert_{\Gamma} = 0$ where $\Gamma$ is a part of $\partial \Omega$, see e.g., \cite{KlibanovNguyenTran:JCP2022, NguyenLiKlibanov:2019}. These Carleman estimates can be used to solve quasilinear elliptic PDEs given the data on only a part of $\partial \Omega$.

We are now in the position to establish the Carleman contraction principle for Problem \ref{p}. We temporarily consider the case when the nonlinearity $F$ is Lipschitz continuous with respect to the second and third variables; i.e., there is a constant $C_F$ depending on $F$ such that
\begin{equation}
    |F(\x, s_1, \p_1) - F(\x, s_2, \p_2)| \leq C_F (|s_1 - s_2| + |\p_1 - \p_2|)
    \label{8}
\end{equation}
for all $s_1, s_2 \in \R$ and $\p_1, \p_2 \in \R^d.$
The Lipschitz continuity will be relaxed later by using a truncation technique; see Remark \ref{remNoLip}.
Let $p > \lceil d/2 \rceil + 2$ be an integer such that $H^p(\Omega)$ can be continuously embedded into $C^2(\overline \Omega).$
We assume that the true solution $u^*$ of \eqref{main_eqn_exact} belongs to $H^p(\Omega).$
Fix $\beta \geq \beta_0$ and $\lambda_0$ as in Lemma \ref{lem:Carleman} such that  Carleman estimate \eqref{Car est} holds true for all $\lambda > \lambda_0$.
Let $H$ be the set of admissible solutions
\begin{equation}
    H = \big\{
        \varphi \in H^p(\Omega): \varphi|_{\partial \Omega} = f, 
        A\nabla \varphi\cdot \nu|_{\partial \Omega} = g
    \big\}.
    \label{H}
\end{equation}
Assume that $H \not = \emptyset.$ 
Define $\Phi: H \to H$ as
\[
    \Phi(u) = \underset{\varphi \in H}{\rm argmin} \, J_u(\varphi)
\]
where
\begin{equation}
    J_u(\varphi) =  \int_{\Omega} 
     e^{2\lambda \mu_\beta(\x)}\big\vert {\rm Div} (A(\x) \nabla \varphi(\x)) + F(\x, u(\x), \nabla u(\x))\big\vert^2d\x
     + \epsilon \Vert \varphi\Vert^2_{H^p(\Omega)}
\end{equation}
for all $u \in H^p(\Omega)$
where $\epsilon > 0$ is the regularization parameter.
\begin{remark}[The well-definedness of $\Phi$]
    It is not hard to verify that the functional $J_u$ has a unique minimizer $\Phi(u) \in H$ for all function $u \in H.$
    Using the compact embedding theorem from $H^p(\Omega)$ to $H^2(\Omega)$, together with the trace theory, one can verify that
       $H$ is weakly closed in $H^p(\Omega)$ and
       $J_u$ is weakly lower semicontinuous on $H$.
    The presence of the regularization term implies that $J_{\varphi}$ is coercive in the sense that $\lim_{\|\varphi\| \to \infty} J_u(\varphi) = \infty$.
    Therefore, by a standard argument in analysis, we can conclude that $J_{u}$ has a minimizer. The uniqueness of the minimizer is due to the strict convexity of $J_u$.
    \label{rem2}
\end{remark}

In practice, given $u \in H$, we solve the linear least square problem to compute $\Phi(u)$. 
This is because
 the map $\varphi \mapsto {\rm Div} (A(\x) \nabla \varphi(\x)) + F(\x, u(\x), \nabla u(\x))$ is affine with respect to $\varphi$. 
 We can use many packages for this purpose.
In computation, we use the optimization package with the command ``lsqlin" of Matlab to minimize $J_u$ and then obtain $\Phi(u)$.

\begin{remark}[The Carleman quasi-reversibility method]
    Fix a function $u \in H$. Let $\varphi = \Phi(u)$.
    Since $\varphi$ is in $H$ and it minimizes $J_u$, roughly speaking, the function $\varphi$ ``almost" solved.
    \begin{equation}
    \left\{
        \begin{array}{ll}
             {\rm Div} (A(\x) \nabla \varphi(\x) )+ F(\x, u(\x), \nabla u(\x)) = 0 &\x \in \Omega,  \\
             \varphi(\x) = f(\x) &\x \in \partial \Omega,\\
             A(\x)\nabla \varphi(\x) \cdot \nu(\x) = g(\x) &\x \in \partial \Omega
        \end{array}
    \right.
    \label{main_eqn1}
\end{equation}
    Due to the presence of the regularization term $\epsilon \|\varphi\|_{H^p(\Omega)}^2$, we call $\varphi$ the regularized solution to \eqref{main_eqn1}. 
    The method to compute the regularized solution to the linear equation \eqref{main_eqn1} by minimizing $J_u$ is named the Carleman quasi-reversibility method. This name is suggested by the presence of the Carleman weight function in the formula of $J_u$ and by the quasi-reversibility method to solve linear PDEs with Cauchy data. See \cite{LattesLions:e1969} for the original work on the quasi-reversibility method.
    \label{rem_quasi}
\end{remark}

For $\epsilon > 0$, $\beta > \beta_0$ and $\lambda > \lambda_0$, define the norm
\begin{equation}
    \|\varphi\|_{\epsilon, \beta, \lambda} = 
    \Big(\int_{\Omega} e^{2\lambda \mu_{\beta}}|\varphi|^2 + |\nabla \varphi|^2 d\x
    \Big)^{1/2} + \frac{\epsilon}{\lambda} \|\varphi\|_{H^p(\Omega)}.
\end{equation}

The contraction behavior of $\Phi$ is confirmed by the following theorem and its consequence mentioned in Corollary \ref{cor1}.
\begin{Theorem}
    There is a number $C$ depending only on  $\x_0,$ $\Omega,$ $\Lambda$, $\beta$ and $d$ such that
    \begin{equation}
        \Vert\Phi(u) - \Phi(v)\Vert_{\epsilon, \beta, \lambda} \leq \sqrt{\frac{C_{F}}{C\lambda}} \Vert u - v\Vert_{\epsilon, \beta, \lambda}
        \label{contraction}
    \end{equation}
    for all $u, v \in H^p(\Omega).$
\label{thm1}
\end{Theorem}

\begin{Corollary}
    Choose $\lambda \gg 1$ such that $\theta = \sqrt{\frac{ C_{F}}{C \lambda}} \in (0, 1)$. It follows from  \eqref{contraction} $\Phi$ is a contraction map with respect to the norm $\|\cdot\|_{\epsilon, \beta, \lambda}$.
    \label{cor1}
\end{Corollary}

\begin{proof}[Proof of Theorem \ref{thm1}]
    Define the set of test functions
    \begin{equation}
        H_0 = \big\{
        \varphi \in H^p(\Omega): \varphi|_{\partial \Omega} = 0, 
        A\nabla \varphi\cdot \nu|_{\partial \Omega} = 0
    \big\}.
    \label{H0}
    \end{equation}
    Recall the admissible set of solutions $H$ defined in \eqref{H}.
    Take two arbitrary functions  $u$ and $v$ in $H$.
    Let $u_1 = \Phi(u)$ and $v_1 = \Phi(v)$.
    Since $u_1$ is the minimizer of $J_{u}$ in $H$,  by the variational principle, we have  for all $h \in H_0$
    \begin{multline}
        \Big\langle
            e^{2\lambda \mu_\beta(\x)} \big[\Div (A(\x) \nabla u_1(\x)) + F(\x, u(\x), \nabla u(\x))\big],
            \Div (A(\x)\nabla h(\x))
        \Big\rangle_{L^2(\Omega)}
        \\
        + \epsilon \big\langle
            u_1(\x), h(\x)
        \big\rangle_{H^p(\Omega)}
        = 0.
        \label{10}
    \end{multline}
    Similarly, for all $h \in H_0$,
    \begin{multline}
        \Big\langle
            e^{2\lambda \mu_\beta(\x)} \big[\Div (A(\x) \nabla v_1(\x)) + F(\x, v(\x), \nabla v(\x))\big],
            \Div (A(\x)\nabla h(\x))
        \Big\rangle_{L^2(\Omega)}
               \\
        + \epsilon \big\langle
            v_1(\x), h(\x)
        \big\rangle_{H^p(\Omega)}
        = 0.
    \label{11}
    \end{multline}
    Combining \eqref{10} and \eqref{11}, using the inequality $2ab \leq a^2 + b^2$ and taking the test function 
    \[
        h = u_1 - v_1 \in H_0,
    \] we have        
        \begin{multline}
        \frac{1}{2} \int_{\Omega} e^{2\lambda \mu_{\beta}} |\Div(A(\x)\nabla h(\x))|^2d\x
        + \epsilon \|h\|_{H^p(\Omega)}^2
        \\ 
        \leq 
        \frac{1}{2} \int_{\Omega} e^{2\lambda \mu_{\beta}} |F(\x, u(\x)  , \nabla u(\x))) - F(\x, v(\x)  , \nabla v(\x)))|^2d\x.
        \label{12}
        \end{multline}
        Using \eqref{8} and \eqref{12}, we have
        \begin{multline}
        \frac{1}{2}\int_{\Omega} e^{2\lambda \mu_{\beta}} |\Div(A(\x)\nabla h(\x))|^2d\x
        + \epsilon \|h\|_{H^p(\Omega)}^2
        \\ 
        \leq 
        \frac{1}{2}C_{F}\int_{\Omega} e^{2\lambda \mu_{\beta}} \big(|u(\x) - v(\x)| + |\nabla u(\x) - \nabla v(\x)|\big)^2d\x.
        \label{13}
        \end{multline}        
        Note that $h|_{\partial \Omega} = 0$ and $A\nabla h \cdot \nu|_{\partial \Omega} = 0.$ We apply the Carleman estimate \eqref{Car est} for $h$ to get
        \begin{equation}
		\int_{\Omega} e^{2\lambda \mu_\beta(\x)}\vert{\rm Div}(A \nabla h) \vert^2d\x
		\geq
		 C\lambda  \int_{\Omega}  e^{2\lambda \mu_\beta(\x)}\vert\nabla h(\x)\vert^2\,d\x
		+ C\lambda^3  \int_{\Omega}   e^{2\lambda \mu_\beta(\x)}\vert h(\x)\vert^2\,d\x	
		\label{14}	
	\end{equation}
	where $C = C( \x_0,  \Omega, \Lambda, d, \beta) > 0$ depends only on the listed parameters.
	Combining \eqref{13} and \eqref{14}, we have
	\begin{multline}
	    C \lambda  \int_{\Omega}  e^{2\lambda \mu_\beta(\x)}\vert\nabla h(\x)\vert^2\,d\x
		+ C \lambda^3  \int_{\Omega}   e^{2\lambda \mu_\beta(\x)}\vert h(\x)\vert^2\,d\x
		+ \epsilon \|h\|_{H^p(\Omega)}^2
		\\
		\leq C_{F}   \int_{\Omega} e^{2\lambda \mu_{\beta}} \big(|u(\x) - v(\x)| + |\nabla u(\x) - \nabla v(\x)|\big)^2d\x.
	\end{multline}
	Therefore,
    \begin{multline}
	     \int_{\Omega}  e^{2\lambda \mu_\beta(\x)}\vert\nabla u_1(\x) -\nabla v_1(\x)\vert^2\,d\x
		+   \vert u_1(\x) - v_1(\x)\vert^2\,d\x
		+ \frac{\epsilon}{\lambda} \|u_1 - v_1\|_{H^p(\Omega)}^2
		\\
		\leq \frac{C_{F} }{C\lambda} \Big[ \int_{\Omega} e^{2\lambda \mu_{\beta}} \big(|u(\x) - v(\x)|^2 + |\nabla u(\x) - \nabla v(\x)|^2\big)d\x
		+ \frac{\epsilon}{\lambda} \|u - v\|_{H^p(\Omega)}^2\Big].
	\end{multline}
	We have proved \eqref{contraction}.
\end{proof}

Corollary \ref{cor1} guarantees  that when $\lambda$ is sufficiently large, the ``fixed-point" like sequence $\{u_n\}_{n \geq 0} \subset H$ defined as
\begin{equation}
    \left\{
    \begin{array}{ll}
        u_0 \in H,\\
        u_{n} = \Phi(u_{n - 1}) & n \geq 1,
    \end{array}
    \right.
    \label{sequence}
\end{equation}
converges to a function $\overline u \in H$ with respect to the norm $\|\cdot\|_{\epsilon, \beta, \lambda}.$
A question arises if $\overline u$ approximates $u^*.$
An affirmative answer will be given in the next section.

\section{The convergence of the Carleman contraction method}\label{sec3}

Recall that $f$ and $g$ are the noisy versions of the boundary data $f^*$ and $g^*$, respectively.
Let $\delta > 0$ be the noise level. By noise, we mean that we assume 
\begin{equation}
    \inf \{\|e\|_{H^p(\Omega)}: e \in E\} < \delta
    \label{noise}
\end{equation}
where $E = \{e \in H^p(\Omega): e|_{\partial \Omega} = f - f^*, A\nabla e\cdot \nu|_{\partial \Omega} = g - g^*\}$. Note that $E$ is nonempty because $u_0 - u^* \in E.$
Due to \eqref{noise}, there exists a function $e \in E$ such that 
	\begin{equation}
		\|e\|_{H^p(\Omega)} < 2\delta.
	\label{ee}
	\end{equation}
	By the continuous embedding from $H^p(\Omega)$ to $C^2(\overline \Omega),$
	we have
	\begin{equation}
	    \|e\|_{C^2(\overline \Omega)} \leq C \delta
	    \label{e_delta}
	\end{equation}

\begin{remark}
    The existence of the ``error" function $e$ satisfying \eqref{ee} and \eqref{e_delta} implies that the differences $f - f^*$ and $g - g^*$ are traces of smooth functions on $\partial \Omega$. 
    That means the noise must be smooth, which might not be realistic.
    This smoothness condition is significant for proving the convergence theorem; see Theorem \ref{thm2}.
    In practice, one can smooth out the data by many existing methods, e.g., by using the well-known spline curves or the Tikhonov regularization approach.
    However, we can relax this step in the numerical study.
    That means our method is stronger than what we can prove. 
    In our numerical study, we do not have to smooth out the noisy data. We directly compute the desired numerical solutions to \eqref{main_eqn} from the given raw, noisy data
    \begin{equation}
        f = f^*(1 + \delta {\rm rand})
        \quad 
        \mbox{and }
        \quad
        g = g^*(1 + \delta {\rm rand})
    \end{equation}
    where ${\rm rand}$ is a function taking uniformly distributed random numbers in the range $[-1, 1].$
\end{remark}

We have the theorem.
\begin{Theorem}
    Fix $\beta \geq \beta_0$. Recall $\lambda_0$ as in Lemma \ref{lemma carl}.
    Let $\lambda \geq \lambda_0$  be such that \eqref{Car est} holds true and the number 
     $\theta$ in Corollary \ref{cor1} is in $(0, 1)$. 
    Let $\{u_n\}_{n \geq 1} \subset H$ be the sequence defined in \eqref{sequence}. The following statements hold.
    \begin{enumerate}
        \item The sequence $\{u_n\}_{n \geq 1}$ converges in to a function $\overline u \in H$ with respect to the norm $\|\cdot\|_{\epsilon, \beta, \lambda}.$
        \item Let $u^*$ be the solution to \eqref{main_eqn}.  Then,
        \begin{multline}
            \|\overline u - u^*\|_{\epsilon, \beta, \lambda}^2 \leq
            \frac{C}{\lambda} \Big[
            \int_{\Omega} e^{2\lambda \mu_{\beta}(\x)} \Big[ |\Div (A(\x)\nabla e(\x))|^2
		 +
	        |e(\x)|^2 
	        \\
	        + |\nabla e(\x)|^2 
	    \Big] d\x
	    +  \epsilon \|e\|_{H^p(\Omega)}^2 +  \epsilon \|u^*\|_{H^p(\Omega)}^2
            \Big]
            \label{u_star_conve}
        \end{multline}
    \end{enumerate}
    where $C$ is a positive constant depending only on $M$, $F$, $\x_0$,  $\Omega$, $A$, $\beta$ and $d$.
    \label{thm2}
\end{Theorem}

Estimate \eqref{u_star_conve} is interesting in the sense that it, together with \eqref{ee}, guarantees that $\overline u$ tends to $u^*$ as the noise level $\delta$ and the regularization parameter $\epsilon$ tends to $0$.
If $\epsilon = O(\delta^2)$, the convergence rate is Lipschitz.

\begin{proof}[Proof of Theorem \ref{thm2}]
	The first part of the theorem is well-known. 
	We only prove the second part of the theorem. 
	We employ the notation $H_0$ defined in \eqref{H0}.
	Fix $n \geq 1$, since $u_n$ is the minimizer of $J_{u_{n-1}}$ in $H$, for all $h \in H_0$, 
	\begin{multline}
	 \Big\langle
            e^{2\lambda \mu_\beta(\x)} \big[\Div (A(\x) \nabla u_n(\x)) + F(\x, u_{n-1}(\x), \nabla u_{n-1}(\x))\big],
            \\
            \Div (A(\x)\nabla h(\x))
        \Big\rangle_{L^2(\Omega)}
        + \epsilon \Big\langle
            u_n(\x), h(\x)
        \Big\rangle_{H^p(\Omega)}
        = 0.
        \label{25}
	\end{multline}
Since  $u^*$ satisfies \eqref{main_eqn},
	\begin{multline}
	 \Big\langle
            e^{2\lambda \mu_\beta(\x)} \big[\Div (A(\x) \nabla u^*(\x)) + F(\x,  u^*(\x), \nabla  u^*(\x))\big],
            \Div (A(\x)\nabla h(\x))
        \Big\rangle_{L^2(\Omega)}
        \\
        + \epsilon \Big\langle
            u^*(\x), h(\x)
        \Big\rangle_{H^p(\Omega)}
        = \epsilon \Big\langle
            u^*(\x), h(\x)
        \Big\rangle_{H^p(\Omega)}.
        \label{27}
	\end{multline}	
	Combining \eqref{25} and \eqref{27}, we have
	\begin{multline}
	 \Big\langle
            e^{2\lambda \mu_\beta(\x)} \big[\Div (A(\x) \nabla(u_n(\x) -  u^*(\x))) + F(\x,   u_{n - 1}(\x), \nabla  u_{n - 1}(\x)) 
            \\
            - F(\x,  u^*(\x), \nabla  u^*(\x))\big],
            \Div (A(\x)\nabla h(\x))
        \Big\rangle_{L^2(\Omega)}
        + \epsilon \Big\langle
           u_n(\x) -  u^*(\x), h(\x)
        \Big\rangle_{H^p(\Omega)}
        \\
        = -\epsilon \Big\langle
            u^*(\x), h(\x)
        \Big\rangle_{H^p(\Omega)}.
        \label{28}
	\end{multline}	
	Recall that $e$ is the function satisfying \eqref{ee} and \eqref{e_delta}.
	Using the test function 
	\begin{equation}
	    h_n =  u_n - u^* - e \in H_0, \quad \mbox{or } u_n - u^* = h_n + e
	    \label{He}
	\end{equation}  in \eqref{28}, we have
	\begin{multline}
	 \Big\langle
            e^{2\lambda \mu_\beta(\x)} \big[\Div (A(\x) \nabla(h_n(\x) + e(\x)) + F(\x,   u_{n-1}(\x), \nabla   u_{n-1}(\x)) 
            \\
            - F(\x,  u^*(\x), \nabla  u^*(\x))\big],
            \Div (A(\x)\nabla h_n(\x))
        \Big\rangle_{L^2(\Omega)}
        + \epsilon \Big\langle
            h_n(\x) + e(\x), h_n(\x)
        \Big\rangle_{H^p(\Omega)}
        \\
        = -\epsilon \Big\langle
            u^*(\x), h_n(\x)
        \Big\rangle_{H^p(\Omega)}.
        \label{29}
	\end{multline}	
	It follows from \eqref{29} and the inequality $2|ab| \leq 4a^2 + b^2/4$ that
	\begin{multline}
	    \int_{\Omega} e^{2\lambda \mu_{\beta}(\x)} |\Div (A(\x)\nabla h_n(\x))|^2 d\x 
	    + \epsilon\|h_n\|^2_{H^p(\Omega)}
	    \\
	    \leq 
	     C\int_{\Omega} e^{2\lambda \mu_{\beta}(\x)} |\Div (A(\x)\nabla e(\x))|^2d\x
	     \\
	     + C\int_{\Omega} e^{2\lambda \mu_{\beta}(\x)} 
	     |F(\x,   u_{n - 1}(\x), \nabla   u_{n-1}(\x)) 
            - F(\x,  u^*(\x), \nabla  u^*(\x))\big|^2 d\x
            \\
            + C \epsilon \|e\|_{H^p(\Omega)}^2 + C \epsilon \|u^*\|_{H^p(\Omega)}^2.
            \label{301}
	\end{multline}
	Using \eqref{8}, we can estimate the second integral on the right-hand side of \eqref{301} as
	\begin{multline}
	    \int_{\Omega} e^{2\lambda \mu_{\beta}(\x)} 
	     |F(\x,   u_{n-1}(\x), \nabla   u_{n - 1}(\x)) 
            - F(\x,  u^*(\x), \nabla  u^*(\x))\big|^2 d\x  
            \\
         \leq C  \int_{\Omega} e^{2\lambda \mu_\beta(\x)}\Big[ | u_{n-1}(\x) - u^*(\x)|^2 + |\nabla (u_{n-1}(\x) - u^*(\x))|^2 \Big]d\x.
        \label{31}
	\end{multline}
	Combining \eqref{301} and \eqref{31}, we get
	\begin{multline}
	    \int_{\Omega} e^{2\lambda \mu_{\beta}(\x)} |\Div (A(\x)\nabla h_n(\x))|^2 d\x 
	    + \epsilon\|h_n\|^2_{H^p(\Omega)}
	    \leq 
	     C\int_{\Omega} e^{2\lambda \mu_{\beta}(\x)} |\Div (A(\x)\nabla e(\x))|^2d\x
	     \\
	     +C  \int_{\Omega} e^{2\lambda \mu_\beta(\x)}\Big[ | u_{n-1}(\x) - u^*(\x)|^2 + |\nabla (u_{n-1}(\x) - u^*(\x))|^2 \Big]d\x
            \\
            + C \epsilon \|e\|_{H^p(\Omega)}^2 + C \epsilon \|u^*\|_{H^p(\Omega)}^2.
            \label{30}
	\end{multline}

	Applying the Carleman estimate \eqref{Car est} for the function $h_n$,
	we obtain
	\begin{equation}
		\int_{\Omega} e^{2\lambda \mu_\beta(\x)}\vert{\rm Div}(A \nabla h_n) \vert^2d\x
		\geq
		 C\lambda  \int_{\Omega}  e^{2\lambda \mu_\beta(\x)}\vert\nabla h_n(\x)\vert^2\,d\x
		+ C\lambda^3  \int_{\Omega}   e^{2\lambda \mu_\beta(\x)}\vert h_n(\x)\vert^2\,d\x.
		\label{34}
	\end{equation}
	Combining \eqref{30} and \eqref{34} and recalling $\lambda \gg 1$, we have
	\begin{multline}
	    \lambda\Big[  \int_{\Omega}  e^{2\lambda \mu_\beta(\x)}\vert\nabla h_n(\x)\vert^2\,d\x
		+ \int_{\Omega}   e^{2\lambda \mu_\beta(\x)}\vert h_n(\x)\vert^2\,d\x\Big]
		\leq 
		 C\int_{\Omega} e^{2\lambda \mu_{\beta}(\x)} |\Div (A(\x)\nabla e(\x))|^2d\x
	     \\
	     +C  \int_{\Omega} e^{2\lambda \mu_\beta(\x)}\Big[ | u_{n-1}(\x) - u^*(\x)|^2 + |\nabla (u_{n-1}(\x) - u^*(\x))|^2 \Big]d\x
            \\
            + C \epsilon \|e\|_{H^p(\Omega)}^2 + C \epsilon \|u^*\|_{H^p(\Omega)}^2.
	    \label{36}
	\end{multline}
	Let $n \to \infty$ and recall that $\{u_n\}_{n\geq 0}$ strongly converges to $\overline u$ in $H_{\lambda, \beta}$. 
	Due to \eqref{He}, we have
	\begin{multline}
	    \lambda\Big[  \int_{\Omega}  e^{2\lambda \mu_\beta(\x)}\Big(\vert\nabla (\overline u(\x) - u^*(\x) - e(\x)) \vert^2 
		+   \vert \overline u(\x) - u^*(\x) - e(\x)\vert^2 \Big)d\x\Big]
		\\
		\leq 
		 C\int_{\Omega} e^{2\lambda \mu_{\beta}(\x)} |\Div (A(\x)\nabla e(\x))|^2 
		d\x
	     \\
	     +C  \int_{\Omega} e^{2\lambda \mu_\beta(\x)}\Big[ | \overline u(\x) - u^*(\x)|^2 + |\nabla (\overline u(\x) - u^*(\x))|^2 \Big]d\x
            \\
            + C \epsilon \|e\|_{H^p(\Omega)}^2 + C \epsilon \|u^*\|_{H^p(\Omega)}^2.
	    \label{36p}
	\end{multline}
	Estimate \eqref{u_star_conve} is a direct consequence of   \eqref{36p} and the inequality $(a - b)^2 \geq \frac{1}{2}a^2 - b^2$.
\end{proof}


\section{Numerical study}\label{sec4}

We consider the case $d = 2$ and $A$ the identity matrix for simplicity.
The computational domain $\Omega$ is chosen to be $(-1, 1)^2$.
We solve \eqref{main_eqn} by the finite difference method.
That means we compute the values of the solution $u^*$ on the grid
\begin{equation}
    \mathcal{G} = \big\{
        (x_i = -1 + (i-1)d_{\x}, y_j = -1 + (j - 1)d_{\x}):
        1 \leq i, j \leq N
    \big\}
    \label{gridG}
\end{equation} where $d_\x = \frac{2}{N-1}$ and $N$ is a large integer. In our numerical study, $N = 150.$

Theorem \ref{thm1} and Theorem \ref{thm2} guarantee that $u_n$, see \eqref{sequence} with $n$ sufficiently large, is an approximation of $u^*$.
This suggests a procedure to compute $u^*$.
This procedure is written in Algorithm \ref{alg}.

\begin{algorithm}[h!]
\caption{\label{alg}The procedure to compute the numerical solution to \eqref{main_eqn}}
	\begin{algorithmic}[1]
	\State \label{s1} Choose a regularization parameter $\epsilon$ and a threshold number $\kappa_0 > 0$.
	 \State \label{s2} Set $n = 0$. Choose an arbitrary initial solution $u_0 \in H.$
		\State \label{step update} 		
		Compute $u_{n + 1} = \Phi(u_n)$ by minimizing $J_{u_{n}}$ in $H$   
	\If {$\|u_{n + 1} - u_n\|_{L^2(\Omega)} > \kappa_0$}
		\State Replace $n$ by $n + 1.$
		\State Go back to Step \ref{step update}.
	\Else	
		\State Set the computed solution $u_{\rm comp} = u_{n + 1}.$
	\EndIf
\end{algorithmic}
\end{algorithm}

The numerical scheme in Algorithm \ref{alg} to solve quasi-linear PDEs with Cauchy boundary data was used when we numerically studied a coefficient inverse problem in \cite{Nguyen:CAMWA2020}. In \cite{Nguyen:CAMWA2020}, we only observed the convergence numerically. The rigorous proof of the convergence was missing.
The convergence of this scheme was partly proved in \cite{LeNguyen:jiip2022} and \cite{NguyenNguyenTruong:arxiv2022}. 
By ``partly", we mean that we only prove that the scheme delivers a numerical solution in a small neighborhood of the true solution. 
However, the convergence of the sequence $\{u_n\}_{n \geq 0}$ to a function $\overline u$ is not guaranteed. 
There might be the case when the sequence $\{u_n\}_{n = 0}^{\infty}$ has two subsequences converging to two different functions.
The new point in the current paper is that this is the first time we can define a contraction mapping to guarantee that the divergence above cannot happen.

We manually choose $\epsilon$ and $\kappa_0$ in Step \ref{s1} by a trial and error process. 
That means we take a reference test in which we know the true solution. Then, we choose $\epsilon$ and $\kappa_0$ such that Algorithm \ref{alg} delivers acceptable numerical solution with noiseless data; i.e. $\delta = 0$. Then, we use these parameters for all other tests and noise levels $\delta.$
The reference test is test 1 below. 
In all of our numerical results, $\epsilon = 10^{-6}$ and $\kappa = 10^{-3}$. The Carleman weight function used in this section is $e^{\lambda |\x - \x_0|^{-\beta}} $ with $\lambda = 3$, $\x_0 = (0, 9)$ and $\beta = 10.$ 
\begin{remark}
The parameters are chosen manually as follows. We take a reference test (Test 1 below) in which we already know the true solution. We then vary these parameters so that the computed solution matches the true one.
Then, we use these parameters for all other tests.
In the process of choosing these artificial parameters, we observe that if $e^{\lambda |\x - \x_0|^{-\beta}}$ is too large (for example in the case $\lambda \gg 1$), the computation is not stable. The computer cannot compute the solution since it might treat some large numbers as $\infty.$
In contrast, if we choose $\lambda$ and $\beta$ such that $e^{\lambda |\x - \x_0|^{-\beta}}$ does not ``numerically" blow up, the computed solutions are satisfactory. 
\end{remark}

\begin{remark}
    In Theorems \ref{thm1} and \ref{thm2}, we need to impose the Lipschitz continuity of the nonlinearity $F$.
    This assumption is essential in proving the convergence of the sequence $\{u_n\}_{n \geq 0}$.
    When $F$ is not Lipschitz, we can employ the truncation technique. 
     Assume that we know in advance that the true solution $u^*$ satisfying $\|u^*\|_{C^1(\overline \Omega)} < M$ for some large number $M$. 
        Define
        \begin{equation}
    \chi_M(\x, s, p)
    = \left\{
        \begin{array}{ll}
             1 & (|s|^2 + |p|^2)^{1/2} \leq M   \\
             \in (0, 1)&  M < (|s|^2 + |p|^2)^{1/2} < 2M\\
             0 & (|s|^2 + |p|^2)^{1/2} \geq  2M
        \end{array}
    \right.
    \label{cutoff}
\end{equation}
and $F_M = \chi_M F.$
        It is obvious that $u^*$ satisfies the problem  
          \begin{equation}
    \left\{
        \begin{array}{ll}
             {\rm Div} (A(\x) \nabla u(\x)) + F_M(\x, u(\x), \nabla u(\x)) = 0 &\x \in \Omega,  \\
             u(\x) = f^*(\x) &\x \in \partial \Omega,\\
             A(\x)\nabla u(\x) \cdot \nu(\x) = g^*(\x) &\x \in \partial \Omega
        \end{array}
    \right.
    \label{main_eqn2}
\end{equation}
Then, we can compute $u^*$ using Algorithm \ref{alg} for \eqref{main_eqn2}.
\label{remNoLip}
\end{remark}

In Step \ref{s2} of Algorithm \ref{alg}, we need to choose a function $u_0$ in $H$. A natural way to compute such a function is to solve the linear problem, obtained by removing from \eqref{main_eqn} the nonlinearity $F$, by the Carleman quasi-reversibility method, see Remark \ref{rem_quasi}. 
We do not present the numerical implementation to solve linear PDEs using the Carleman quasi-reversibility method in this paper. The reader can find the details about this in \cite{LeNguyen:jiip2022, Nguyen:CAMWA2020, Nguyens:jiip2020}. 

In Step \ref{step update}, we minimize $J_{u_n}$ in $H$. Similarly to the discussion in Remark \ref{rem_quasi}, the obtained minimizer $u_{n+1}$ is actually the regularized solution to 
\begin{equation}
    \left\{
        \begin{array}{ll}
             \Delta u_{n + 1}(\x) + F(\x, u_n(\x), \nabla u_n(\x)) = 0 &\x \in \Omega,  \\
             u_{n + 1}(\x) = f(\x) &\x \in \partial \Omega,\\
             \partial_{\nu} u_{n+1}(\x)  = g(\x) &\x \in \partial \Omega.
        \end{array}
    \right.
    \label{main_eqn40}
\end{equation}
The details in implementation to compute the regularized solution $u_{n + 1}$ were presented in \cite{LeNguyen:jiip2022, Nguyen:CAMWA2020, Nguyens:jiip2020}, in which we employ the optimization package already built in Matlab. We do not repeat it here.
We next display our numerical examples.

\begin{figure}[h!]
    \centering
    \subfloat[The true solution $u^*$  ]{\includegraphics[width=.3\textwidth]{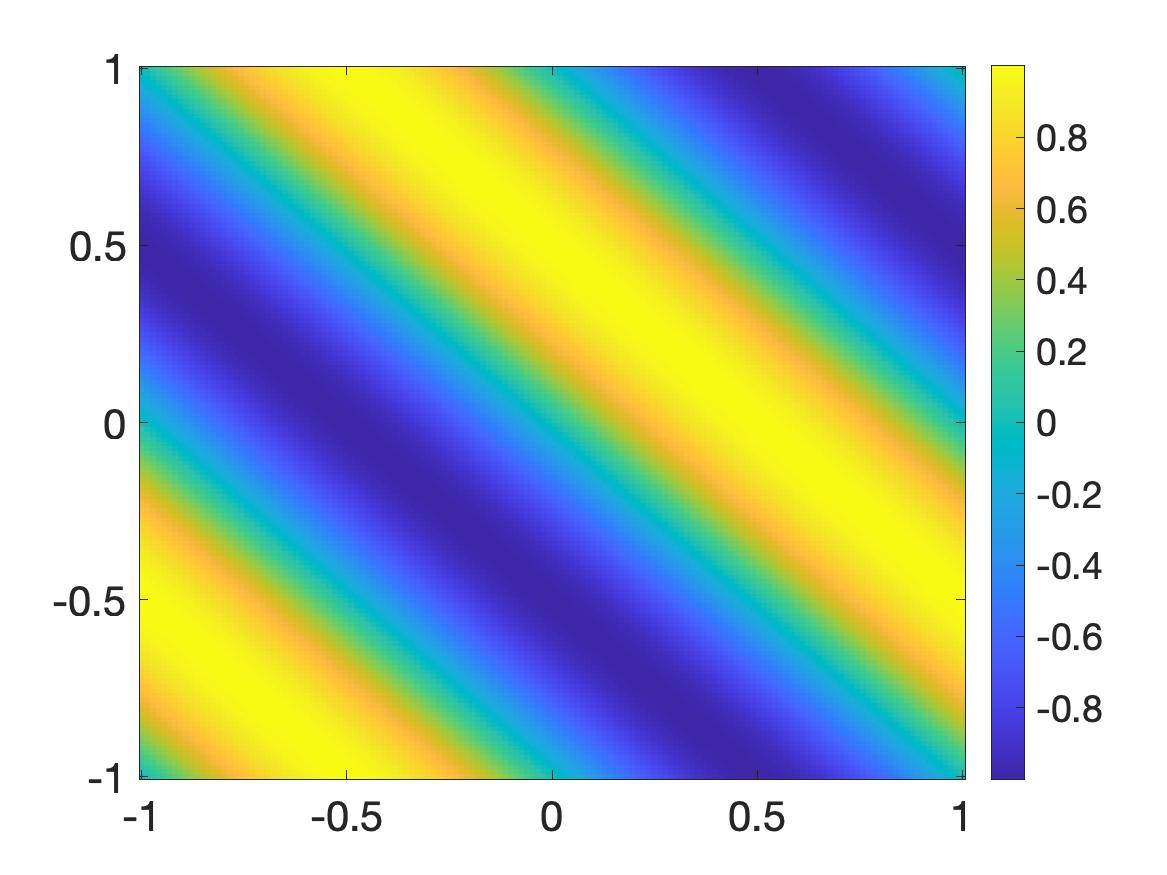}}
    \hfill
    \subfloat[The computed solution $u$  when $\delta = 0\%$ ]{\includegraphics[width=.3\textwidth]{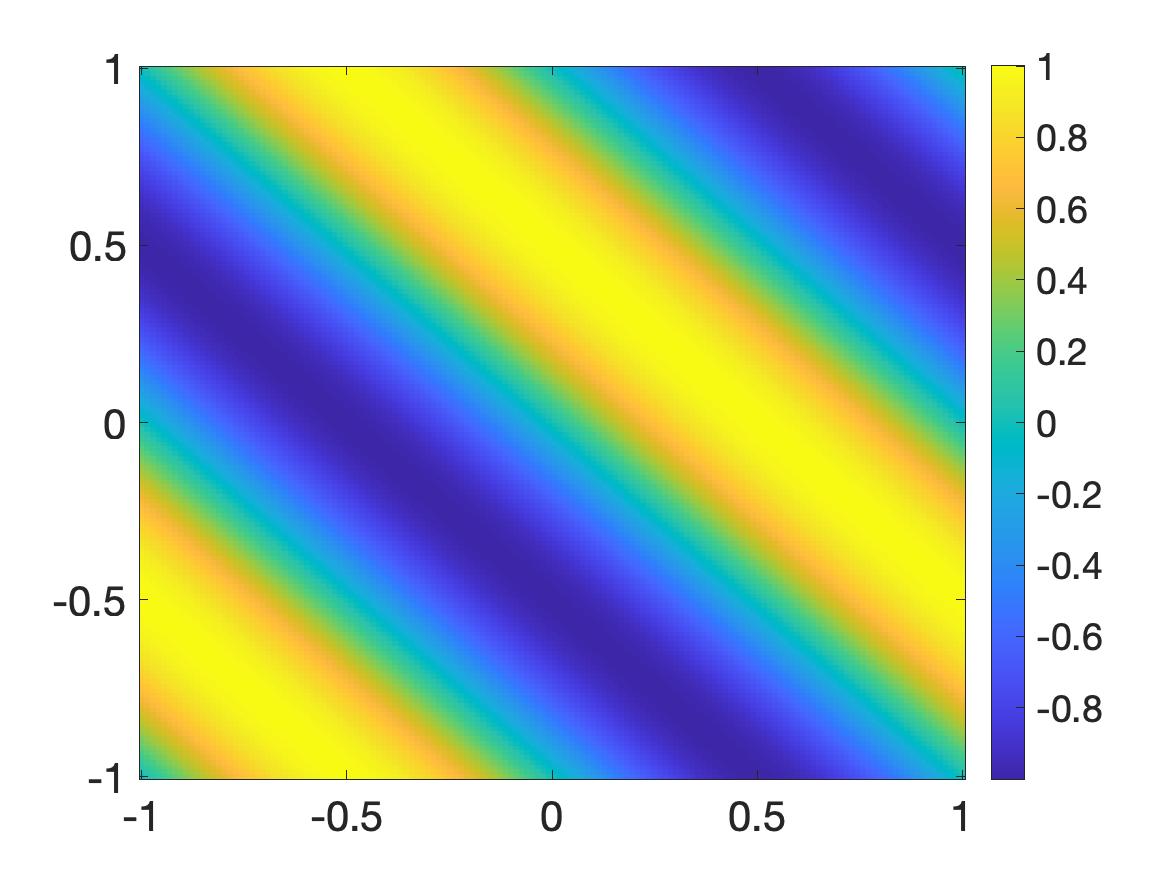}}
    \hfill
    \subfloat[The relative error $\frac{|u^* - u|}{\|u^*\|_{L^\infty(\Omega)}}$ when $\delta = 0\%$ ]{\includegraphics[width=.3\textwidth]{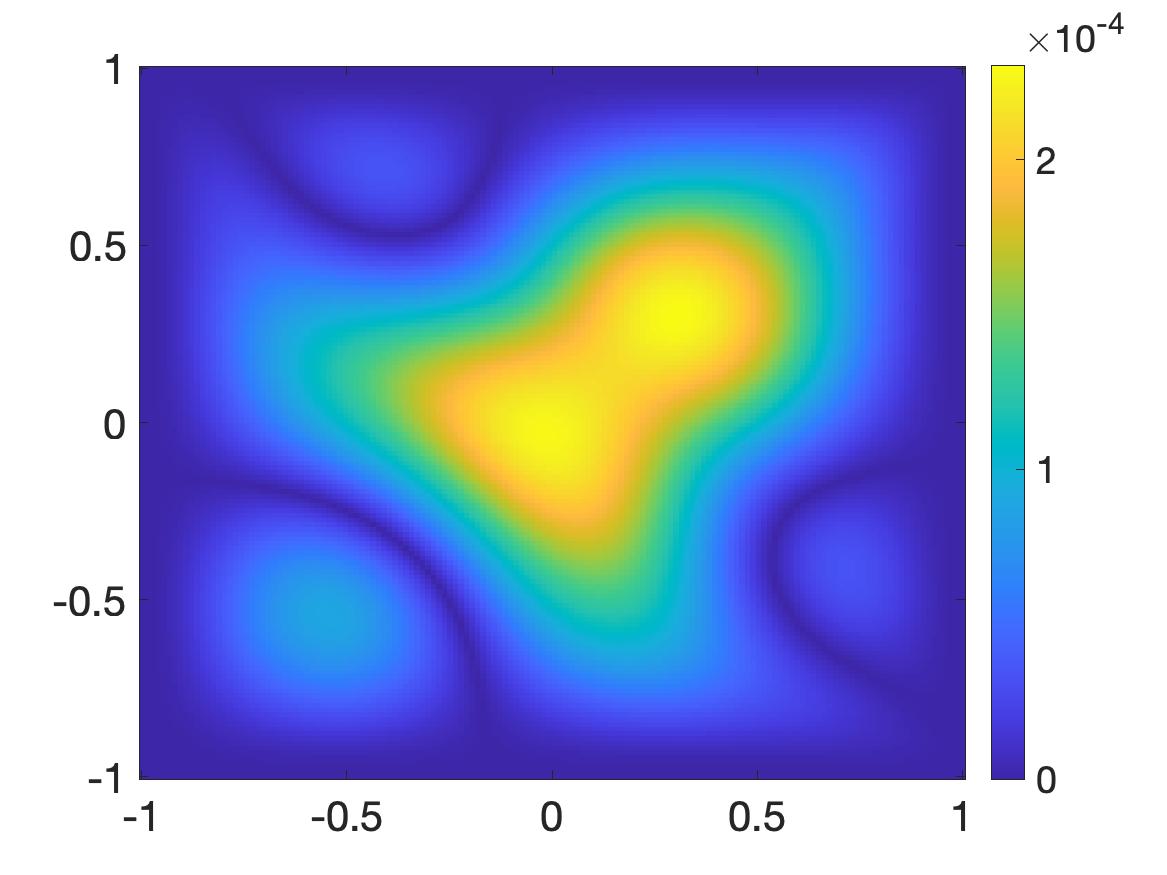}}
    
    \subfloat[The computed solution $u$  when $\delta = 10\%$ ]{\includegraphics[width=.3\textwidth]{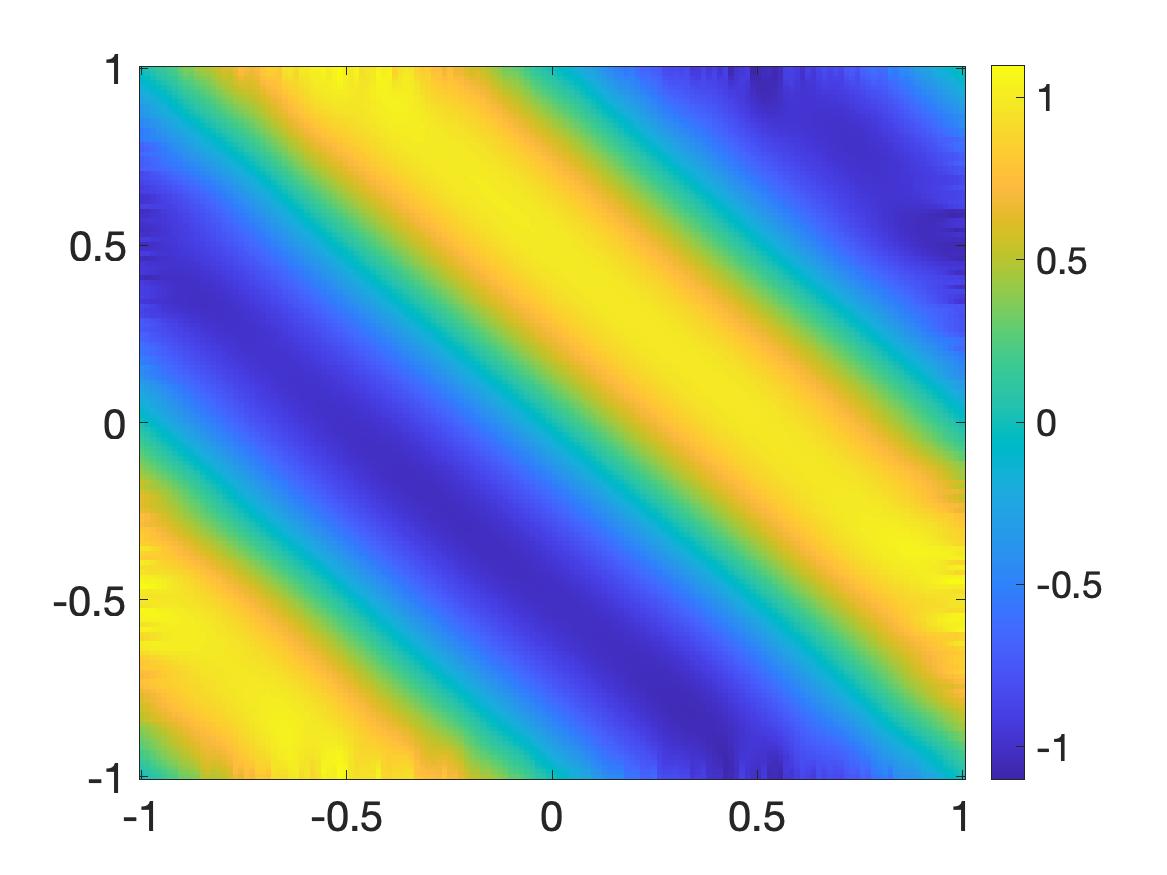}}
    \hfill
    \subfloat[\label{1e}The relative error $\frac{|u^* - u|}{\|u^*\|_{L^\infty(\Omega)}}$ when $\delta = 10\%$ ]{\includegraphics[width=.3\textwidth]{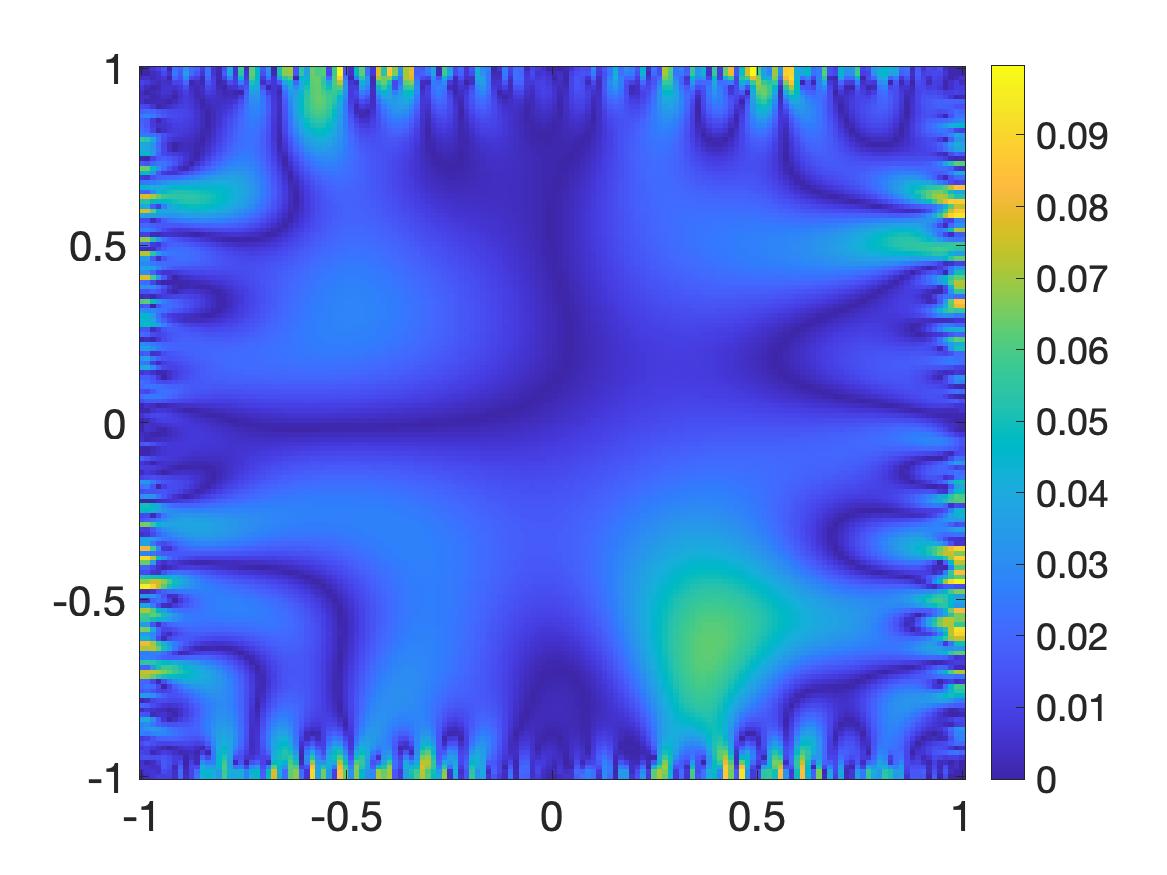}}
    \hfill
    \subfloat[\label{fig1f}The difference $\|u_{n+1} - u_n\|_{L^\infty(\Omega)}$ ]{\includegraphics[width=.3\textwidth]{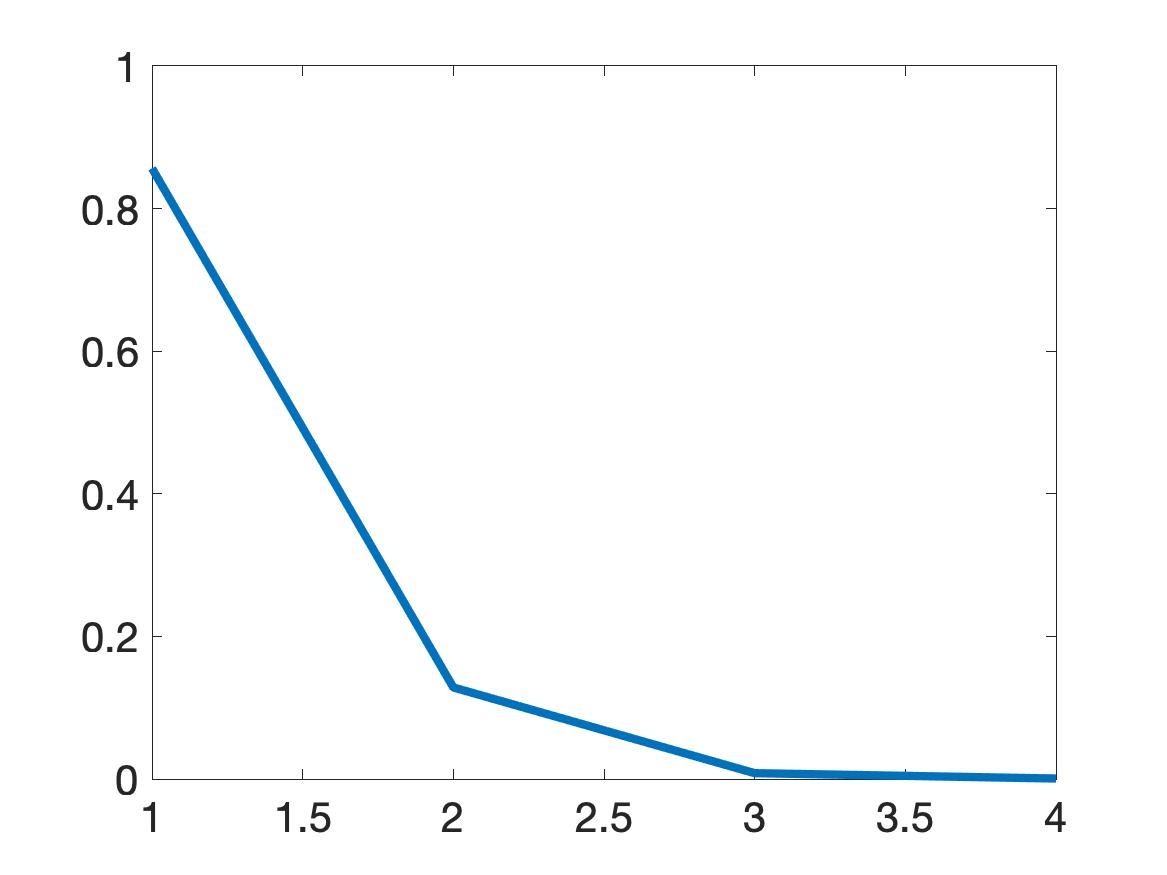}}
    \caption{Test 1. The graphs of the true and computed solution to \eqref{eqn1}, \eqref{f1} and \eqref{g1}  with noiseless and noisy boundary data. }
    \label{fig1}
\end{figure}

\noindent {\bf Test 1. }
In this test, we compute the solution to
\begin{multline}
    \Delta u(\x) + u(\x) + \sqrt{|\nabla u|^2 + 1}
    - \big[
      -2\pi^2 \sin(\pi(x + y)) 
      \\
      + \sin(\pi(x + y))
      +
      \sqrt{\pi^2 \cos(\pi(x + y)) + 1} 
    \big] = 0
    \label{eqn1}
\end{multline}
for all $\x = (x, y) \in \Omega.$
The boundary data are given by
\begin{align}
    u(\x) &= \sin(\pi(x + y))(1 + \delta \mbox{rand}_1), \label{f1}\\
    \partial_{\nu} u(\x)  &= \pi(\cos(\pi(x + y)), \cos(\pi(x + y)))\cdot \nu (1 + \delta \mbox{rand}_2) \label{g1}
\end{align}
for all $\x = (x, y) \in \partial \Omega$,
where $\delta > 0$ is the noise level and $\mbox{rand}_{i}$, $i = 1, 2$, is the function taking uniformly distributed random numbers in the rank $[-1, 1]$.
The true solution to \eqref{eqn1}, \eqref{f1} and \eqref{g1} when $\delta = 0$ is $u^*(\x) = \sin(\pi (x + y))$ for all $\x = (x, y) \in \Omega.$
The graphs of the true and computed solution and the relative $L^{\infty}$ error in the computation are displayed in Figure \ref{fig1}.

It is evident from Figure \ref{fig1} that the numerical solutions to \eqref{eqn1}, \eqref{f1} and \eqref{g1} are computed with a good accuracy.
The relative errors in the computation are in Table \ref{table1}. On the other hand, one can observe from Figure \ref{fig1f} that our method converges fast after only four iterations.

\begin{table}[h!]
    \centering
    \begin{tabular}{c c c c}
    \hline
       Noise level  & $\|u^* - u^{\rm comp}\|_{L^\infty(\Omega)}$ & $\|u^* - u^{\rm comp}\|_{L^2(\Omega)}$  \\
      \hline
      $ \delta = 0\%$  
      &$2.3024\times 10^{-4}$ & 
      $1.2581\times 10^{-4}$ 
      \\
      $ \delta = 2\%$  
      &$0.0200$ & 
      $0.0061$ 
      \\
      $ \delta = 5\%$  
      &$0.0491$ & 
      $0.0153$ 
      \\
      $ \delta = 10\%$  
      &$0.0996$ & 
      $0.0331$ 
      \\
      \hline
    \end{tabular}
    \caption{Test 1. The relative errors in computation.}
    \label{table1}
\end{table}

\begin{figure}[h!]
    \centering
    \subfloat[The true solution $u^*$  ]{\includegraphics[width=.3\textwidth]{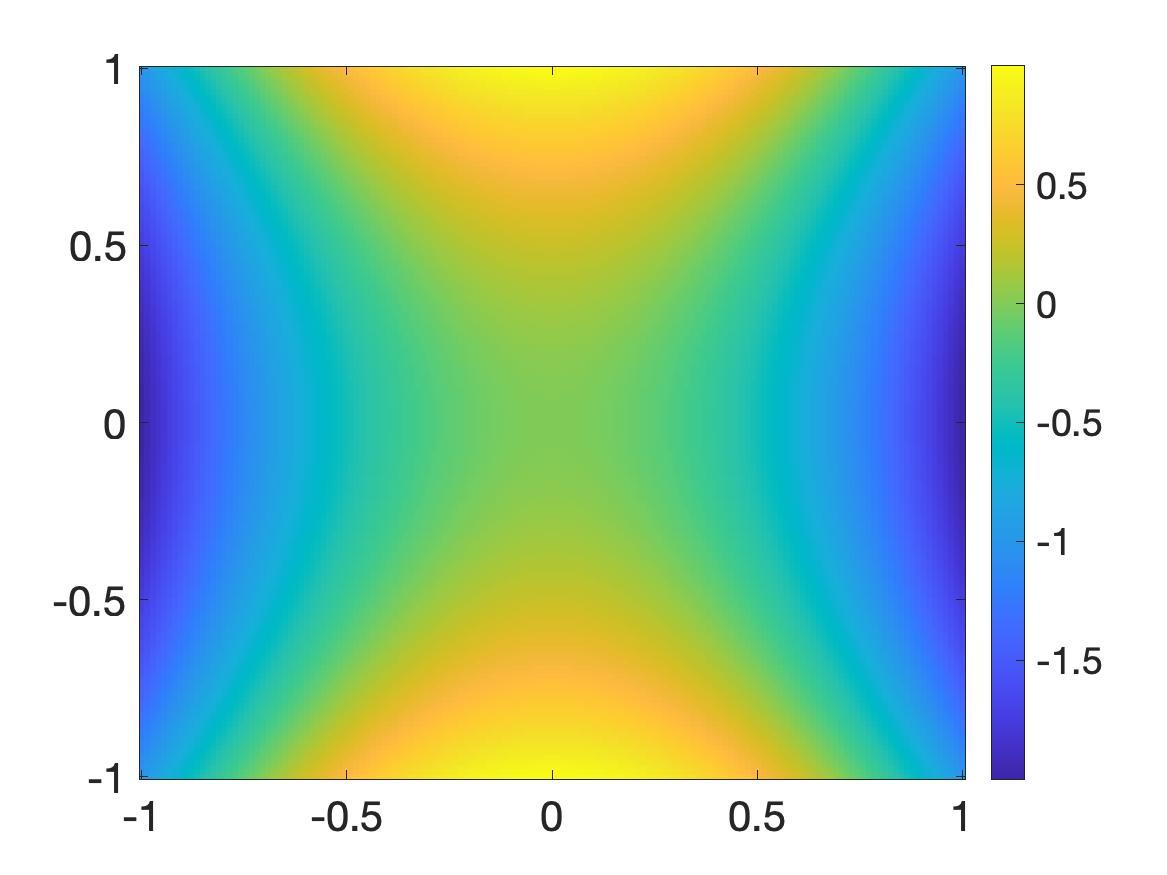}}
    \hfill
    \subfloat[The computed solution $u$  when $\delta = 0\%$ ]{\includegraphics[width=.3\textwidth]{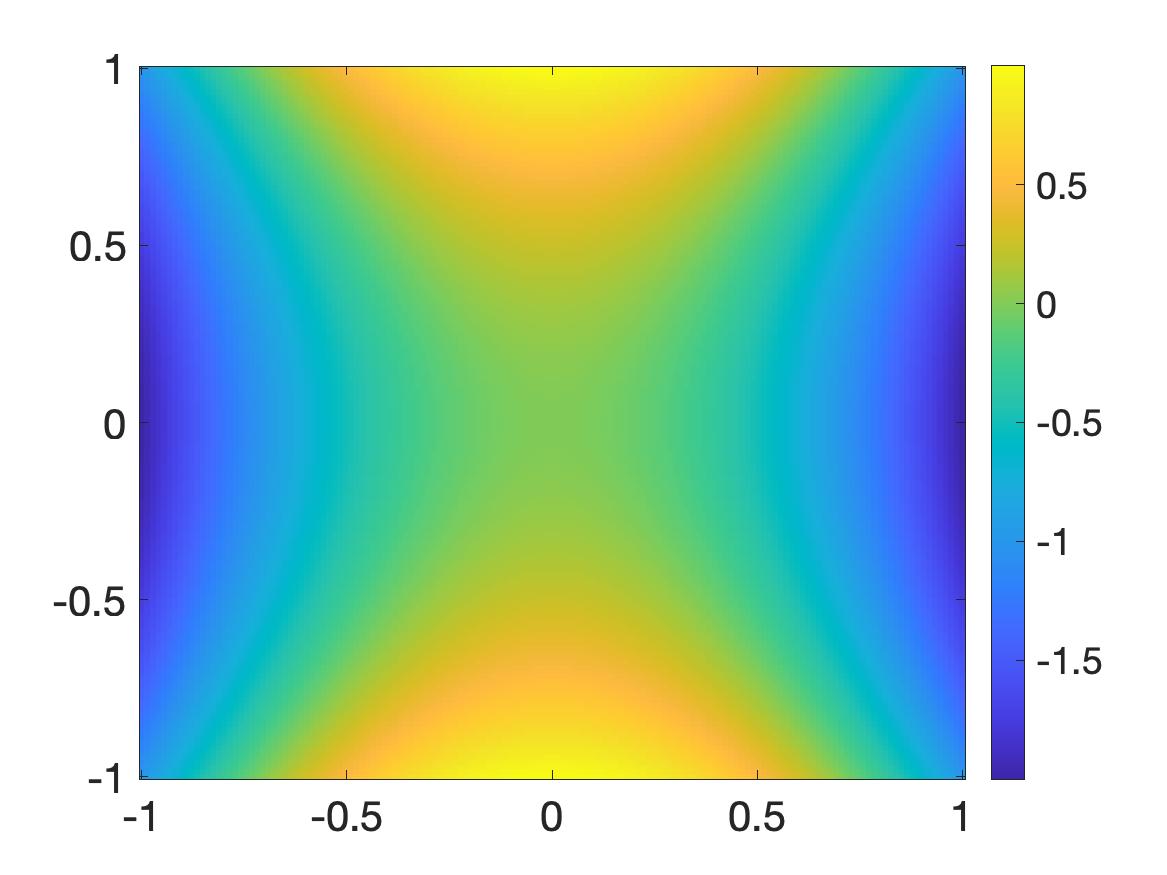}}
    \hfill
    \subfloat[The relative error $\frac{|u^* - u|}{\|u^*\|_{L^\infty(\Omega)}}$ when $\delta = 0\%$ ]{\includegraphics[width=.3\textwidth]{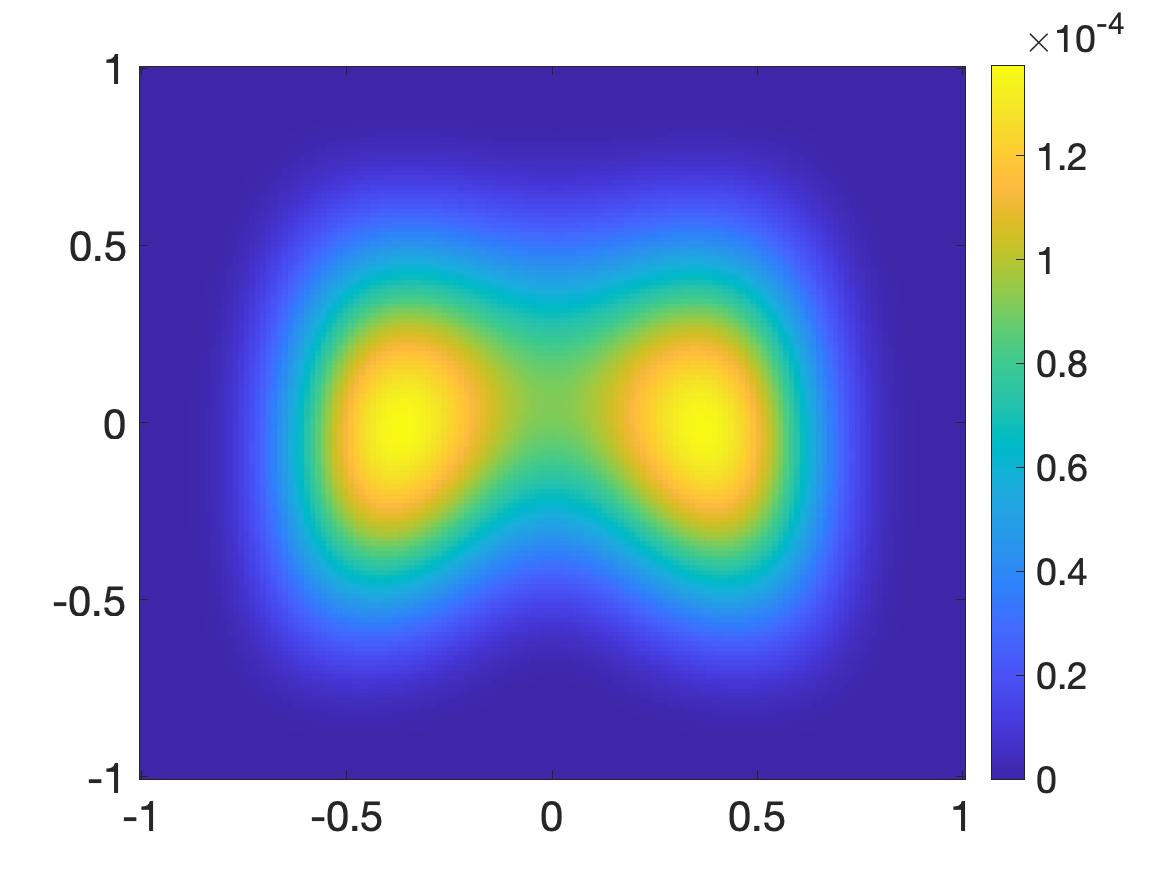}}
    
    \subfloat[The computed solution $u$  when $\delta = 10\%$ ]{\includegraphics[width=.3\textwidth]{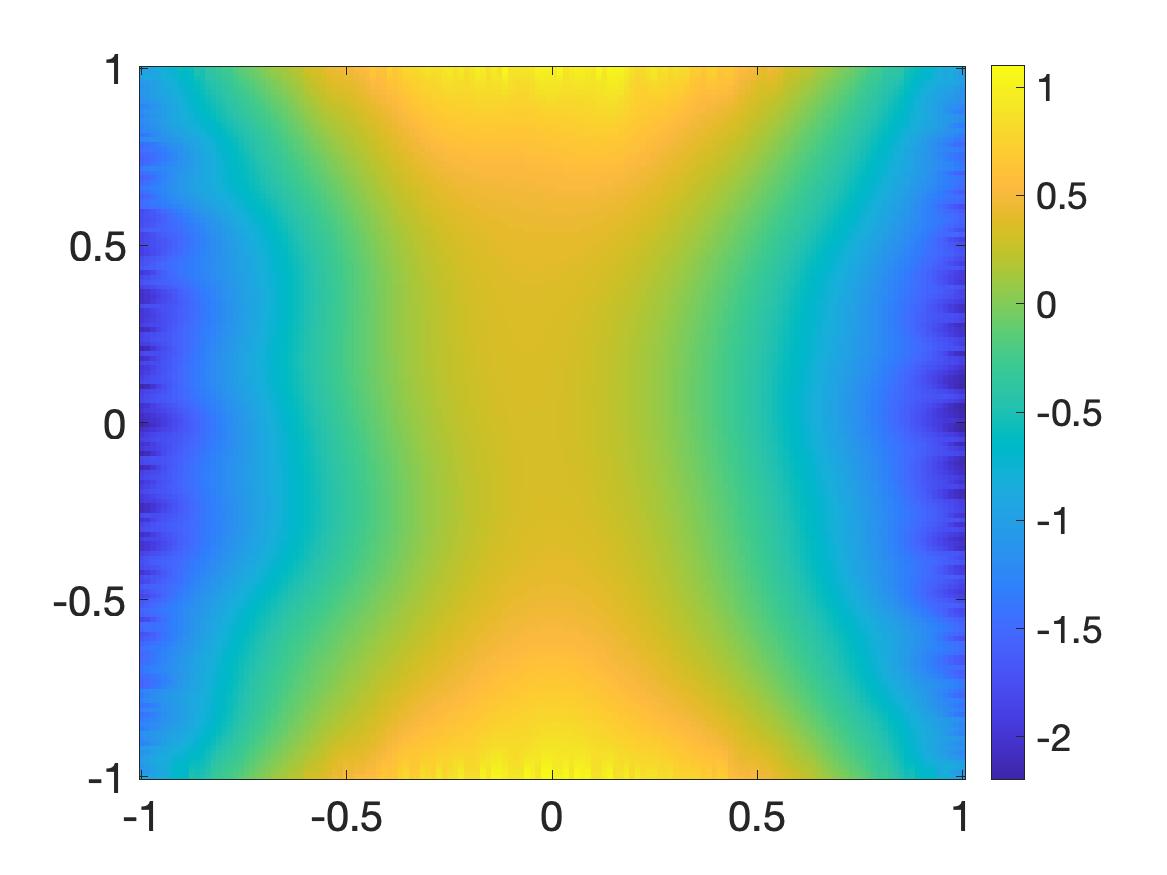}}
    \hfill
    \subfloat[The relative error $\frac{|u^* - u|}{\|u^*\|_{L^\infty(\Omega)}}$ when $\delta = 10\%$ ]{\includegraphics[width=.3\textwidth]{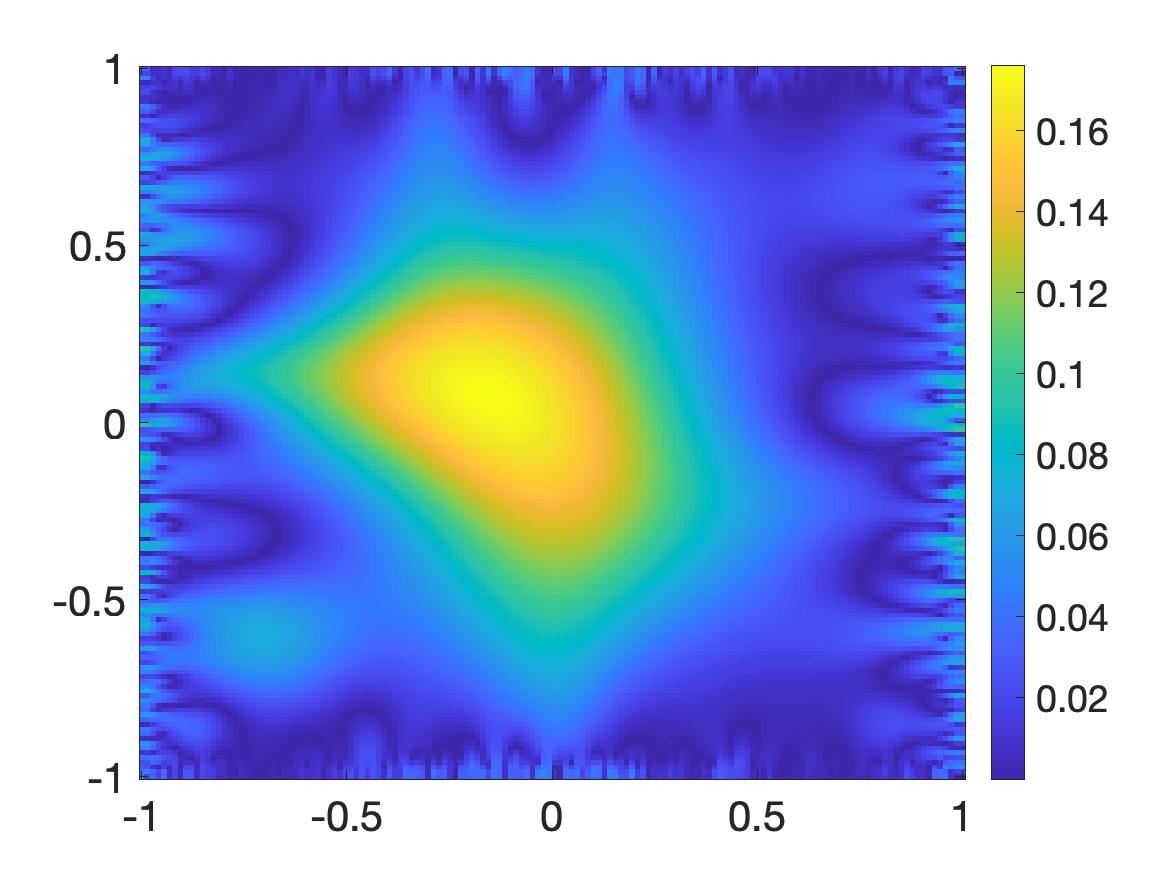}}
    \hfill
    \subfloat[\label{fig2f}The difference $\|u_{n+1} - u_n\|_{L^\infty(\Omega)}$ ]{\includegraphics[width=.3\textwidth]{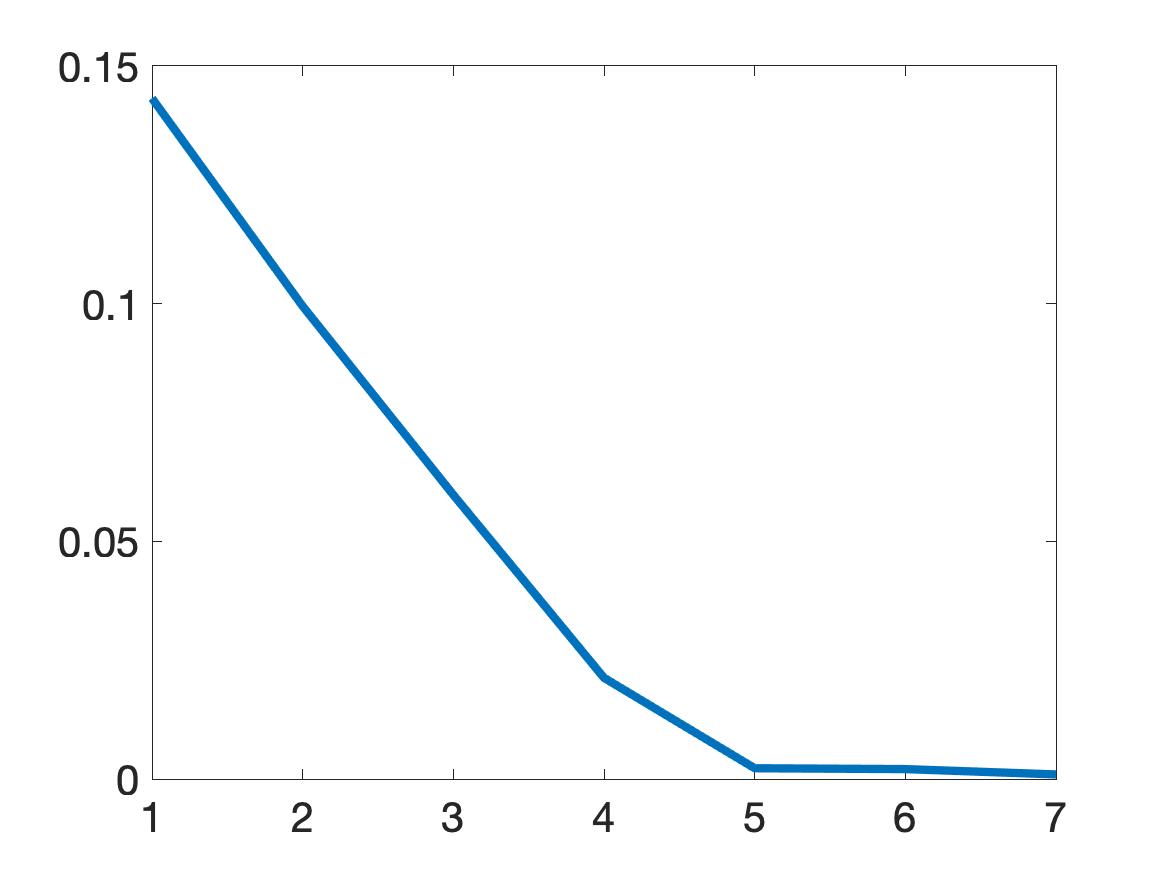}}
    \caption{Test 2. The graphs of the true and computed solution to \eqref{eqn2}, \eqref{f2} and \eqref{g2}  with noiseless and noisy boundary data. }
    \label{fig2}
\end{figure}

\noindent {\bf Test 2. }
We consider a more complicated test with the nonlinearity $F(\x, s, p)$ grows as $O(|p|^2)$ as $p \to \infty$ and is not convex with respect to $p$.
We solve the equation
\begin{equation}
    \Delta u(\x) + u_x - u_y^2  - \big[
       -2 + 2x - 16y^2
    \big] = 0
    \label{eqn2}
\end{equation}
for all $\x = (x, y) \in \Omega.$
The boundary data are given by
\begin{align}
    u(\x) &= (x^2 - 2y^2)(1 + \delta\mbox{rand}_1), \label{f2}\\
    \partial_{\nu} u(\x)  &= (2x, -4y)\cdot \nu (1 + \delta \mbox{rand}_2) \label{g2}
\end{align}
for all $\x = (x, y) \in \partial \Omega$,
where $\delta > 0$ is the noise level and $\mbox{rand}_{i}$, $i = 1, 2$, is the function taking uniformly distributed random numbers in the rank $[-1, 1]$.
The true solution to \eqref{eqn2}, \eqref{f2} and \eqref{g2} when $\delta = 0$ is $u^*(\x) = x^2 - 2y^2$ for all $\x = (x, y) \in \Omega.$
The graphs of the true and computed solution and the relative $L^{\infty}$ error in the computation are displayed in Figure \ref{fig2}.

Even though this test is challenging, Algorithm \ref{alg} delivers out-of-expectation numerical solutions.
The relative errors in the computation are in Table \ref{table2}. On the other hand, one can observe from Figure \ref{fig2f} that our method converges fast after only seven iterations. The number of iterations in this test is greater than that in test 1, probably because the nonlinearity $F$ in this test grows faster at $|p| \to \infty.$ 

\begin{table}[h!]
    \centering
    \begin{tabular}{c c c }
    \hline
       Noise level  & $\|u^* - u^{\rm comp}\|_{L^\infty(\Omega)}$ & $\|u^* - u^{\rm comp}\|_{L^2(\Omega)}$ \\
      \hline
      $ \delta = 0\%$  
      &$1.3727\times 10^{-4}$ & 
      $1.3134\times 10^{-4}$ 
      \\
      $ \delta = 2\%$  
      &$0.0198$ & 
      $ 0.0130$ 
      \\
      $ \delta = 5\%$  
      &$ 0.0702$ & 
      $ 0.0694$ 
      \\
      $ \delta = 10\%$  
      &$0.1760$ & 
      $0.1721$ 
      \\
      \hline
    \end{tabular}
    \caption{Test 2. The relative errors in computation}
    \label{table2}
\end{table}

\begin{figure}[h!]
    \centering
    \subfloat[The true solution $u^*$  ]{\includegraphics[width=.3\textwidth]{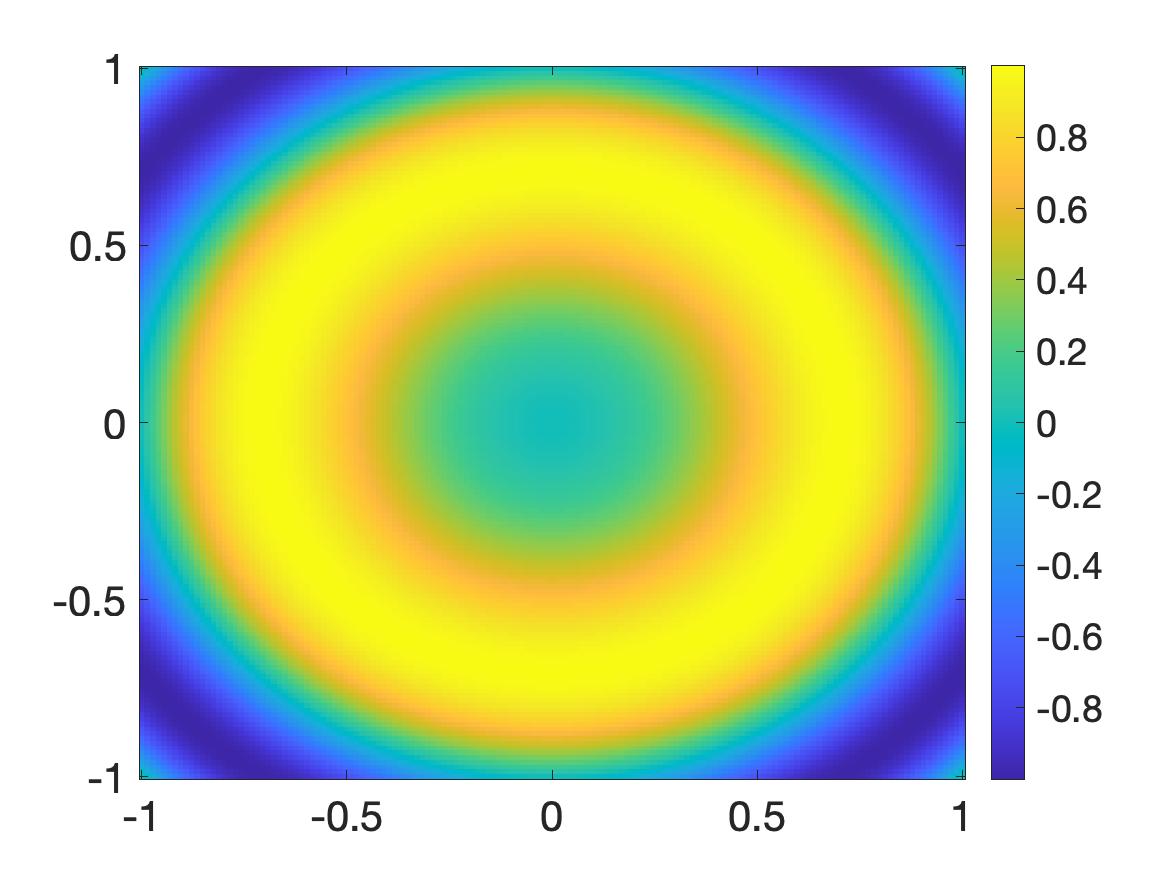}}
    \hfill
    \subfloat[The computed solution $u$  when $\delta = 0\%$ ]{\includegraphics[width=.3\textwidth]{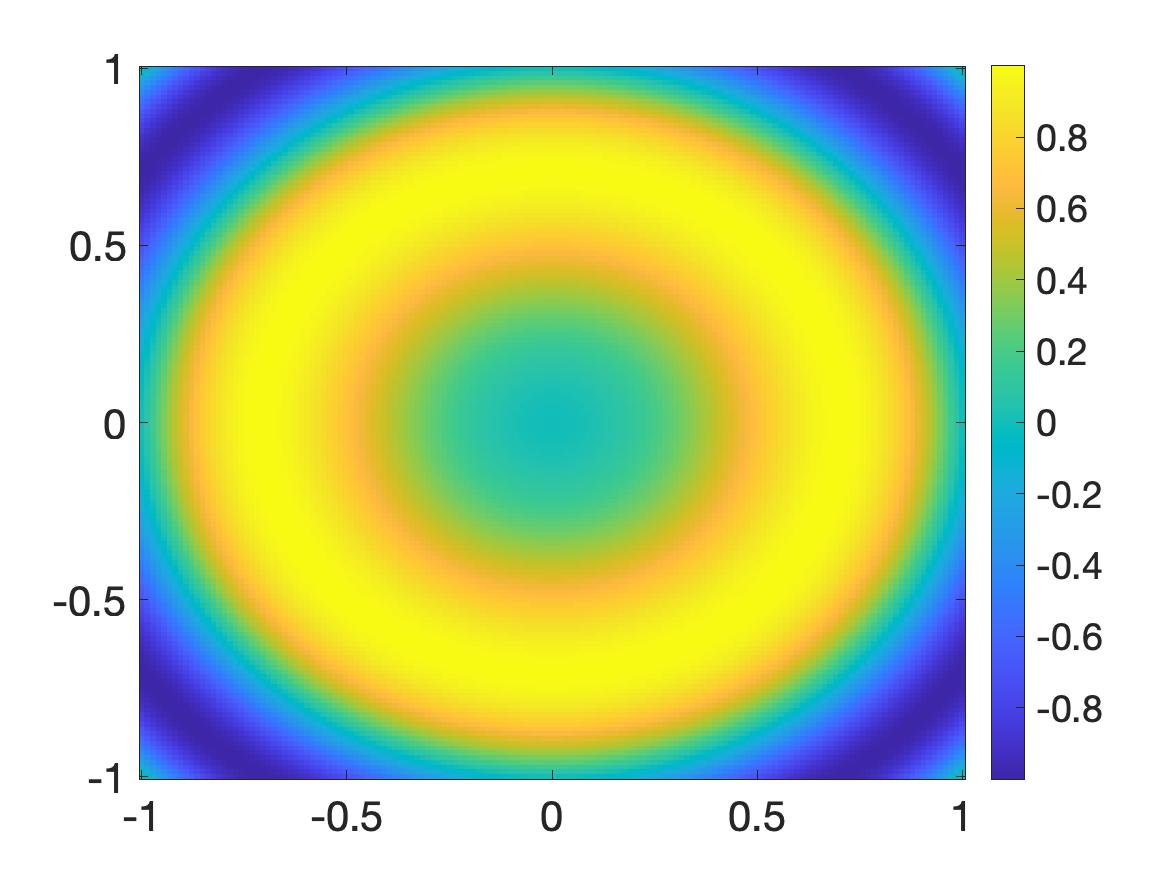}}
    \hfill
    \subfloat[The relative error $\frac{|u^* - u|}{\|u^*\|_{L^\infty(\Omega)}}$ when $\delta = 0\%$ ]{\includegraphics[width=.3\textwidth]{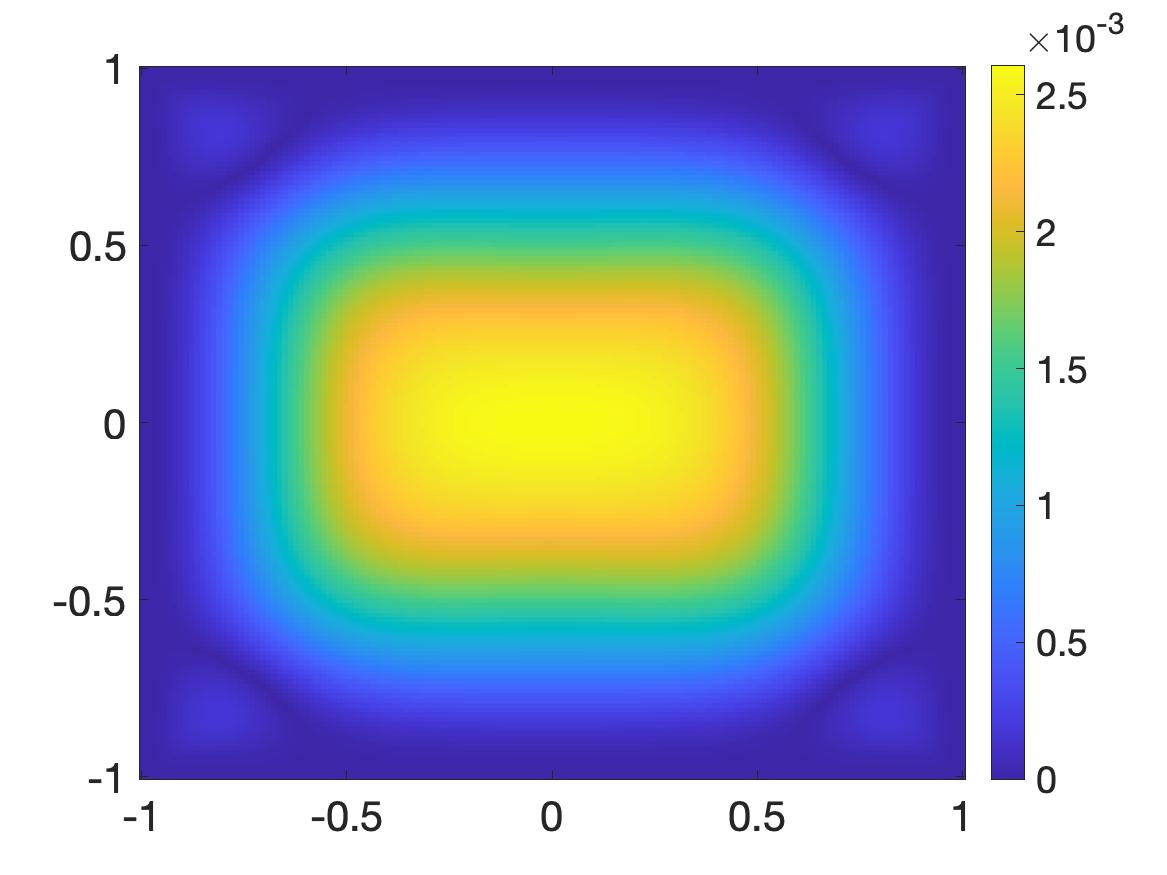}}
    
    \subfloat[The computed solution $u$  when $\delta = 10\%$ ]{\includegraphics[width=.3\textwidth]{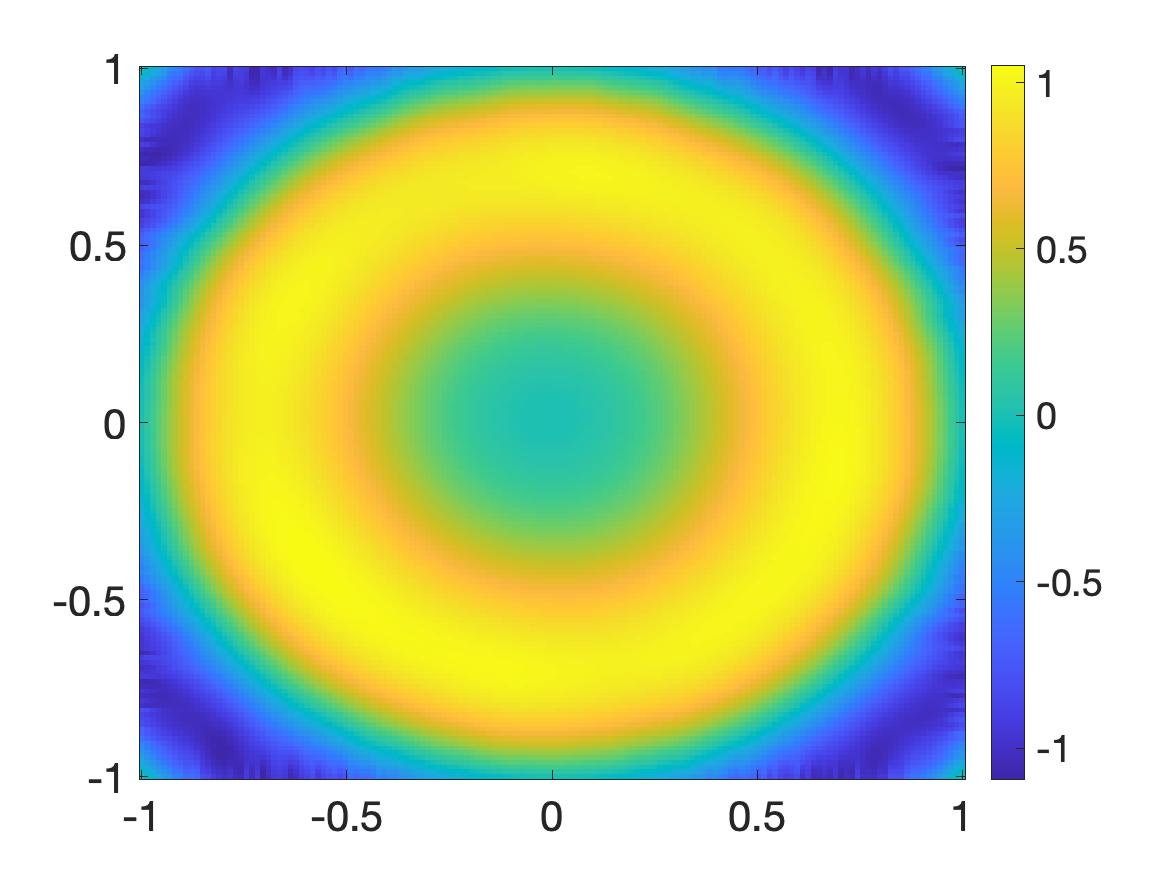}}
    \hfill
    \subfloat[The relative error $\frac{|u^* - u|}{\|u^*\|_{L^\infty(\Omega)}}$ when $\delta = 10\%$ ]{\includegraphics[width=.3\textwidth]{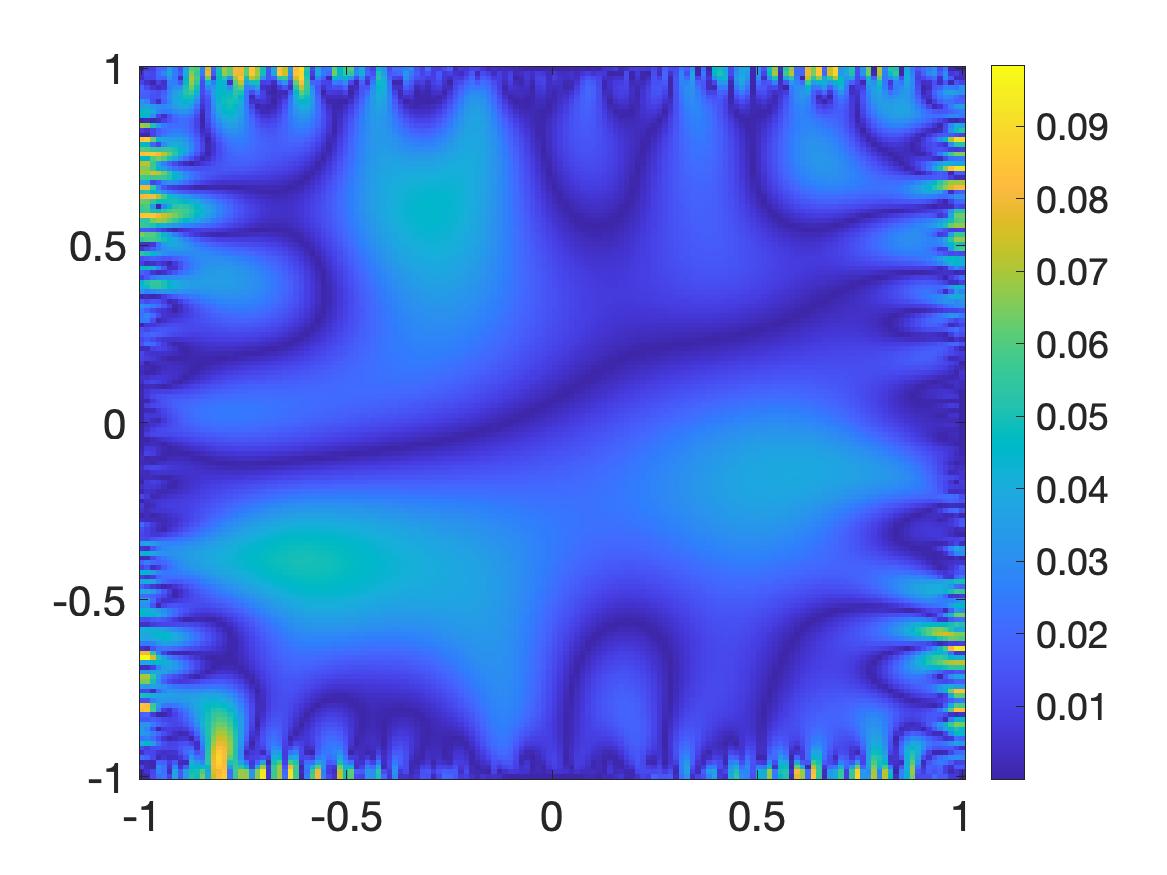}}
    \hfill
    \subfloat[\label{fig3f}The difference $\|u_{n+1} - u_n\|_{L^\infty(\Omega)}$ ]{\includegraphics[width=.3\textwidth]{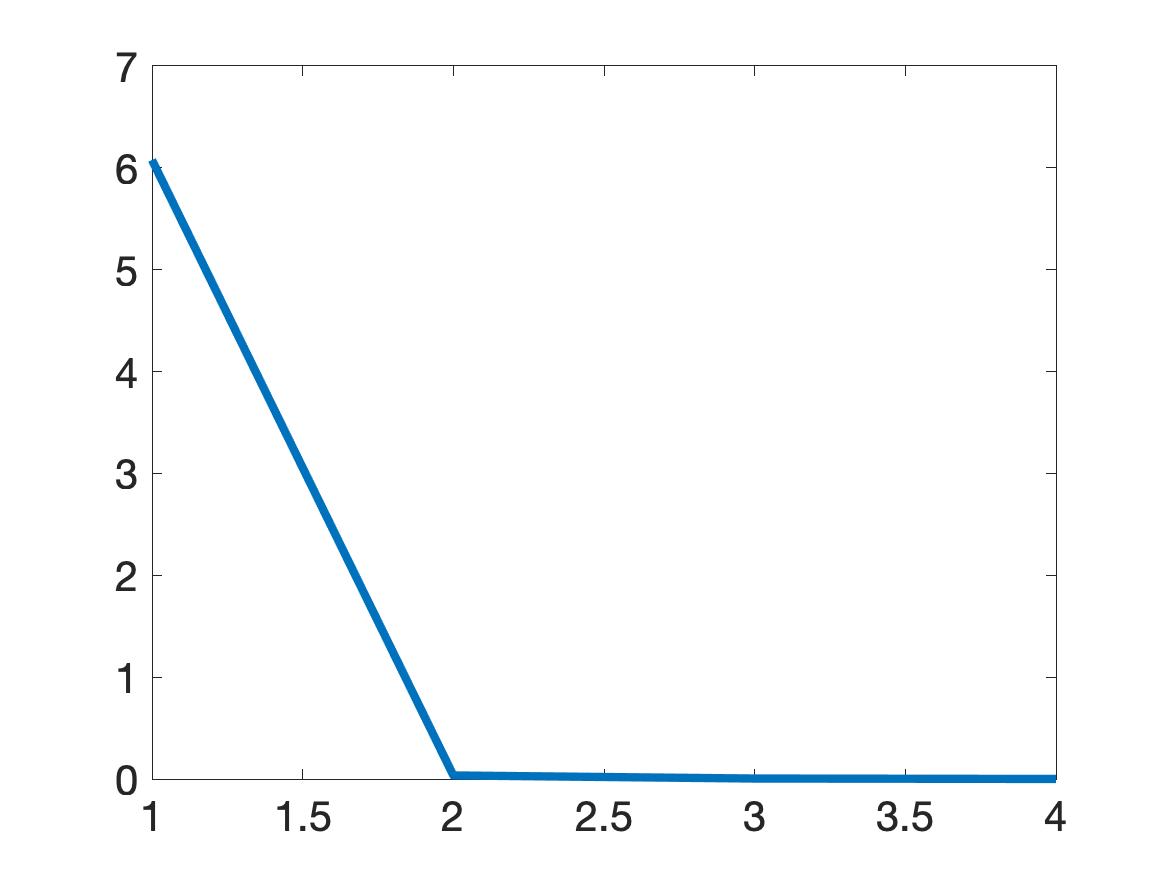}}
    \caption{Test 3. The graphs of the true and computed solution to \eqref{eqn3}, \eqref{f3} and \eqref{g3}  with noiseless and noisy boundary data. }
    \label{fig3}
\end{figure}

\noindent {\bf Test 3. }
In this test, we try the efficiency of the Algorithm when the nonlinearity $F$ is not smooth.
We solve the equation
\begin{multline}
    \Delta u(\x) + |u_x(\x)| - |u_y(\x)|  +
       4\pi\big(
          \pi(x^2 + y^2)\sin(\pi(x^2 + y^2))
          -\cos(\pi(x^2 + y^2))
       \big)
       \\
     - 2\pi\big(|x \cos(\pi(x^2 + y^2))| - |y\cos(\pi(x^2 + y^2))|\big)
     = 0
    \label{eqn3}
\end{multline}
for all $\x = (x, y) \in \Omega.$
The boundary data are given by
\begin{align}
    u(\x) &= \sin(\pi(x^2 + y^2))(1 + \delta\mbox{rand}_1), \label{f3}\\
    \partial_{\nu} u(\x)  &= 2\pi(x\cos(\pi(x^2 + y^2)), y\cos(\pi(x^2 + y^2)))\cdot \nu (1 + \delta \mbox{rand}_2) \label{g3}
\end{align}
for all $\x = (x, y) \in \partial \Omega$,
where $\delta > 0$ is the noise level and $\mbox{rand}_{i}$, $i = 1, 2$, is the function taking uniformly distributed random numbers in the rank $[-1, 1]$.
The true solution to \eqref{eqn3}, \eqref{f3} and \eqref{g3} when $\delta = 0$ is $u^*(\x) = \sin(\pi(x^2 + y^2))$ for all $\x = (x, y) \in \Omega.$
The graphs of the true and computed solution and the relative $L^{\infty}$ error in the computation are displayed in Figure \ref{fig3}.

Even though this test is challenging, Algorithm \ref{alg} delivers out-of-expectation numerical solutions.
The relative errors are compatible with the noise, which can be found in Table \ref{table3}. On the other hand, one can observe from Figure \ref{fig3f} that our method converges fast.
The stopping criterion meets after only four iterations.

\begin{table}[h!]
    \centering
    \begin{tabular}{c c c }
    \hline
       Noise level  & $\|u^* - u^{\rm comp}\|_{L^\infty(\Omega)}$ & $\|u^* - u^{\rm comp}\|_{L^2(\Omega)}$ \\
      \hline
      $ \delta = 0\%$  
      &$0.0026$ & 
      $0.0018$ 
      \\
      $ \delta = 2\%$  
      &$0.0200$ & 
      $ 0.0062$ 
      \\
      $ \delta = 5\%$  
      &$ 0.0509$ & 
      $ 0.0168$ 
      \\
      $ \delta = 10\%$  
      &$0.0983$ & 
      $0.0332$ 
      \\
      \hline
    \end{tabular}
    \caption{Test 3. The relative errors in computation}
    \label{table3}
\end{table}

\begin{figure}[h!]
    \centering
    \subfloat[The true solution $u^*$  ]{\includegraphics[width=.3\textwidth]{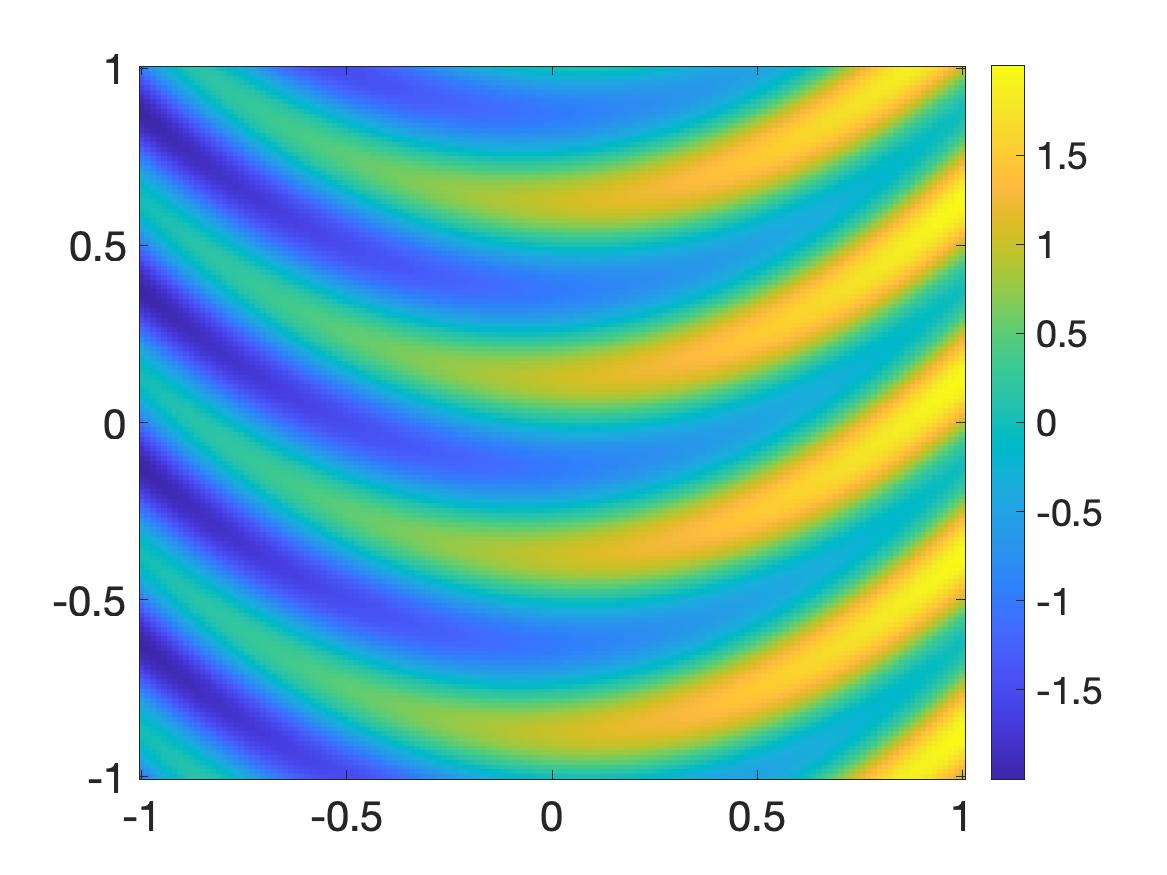}}
    \hfill
    \subfloat[The computed solution $u$  when $\delta = 0\%$ ]{\includegraphics[width=.3\textwidth]{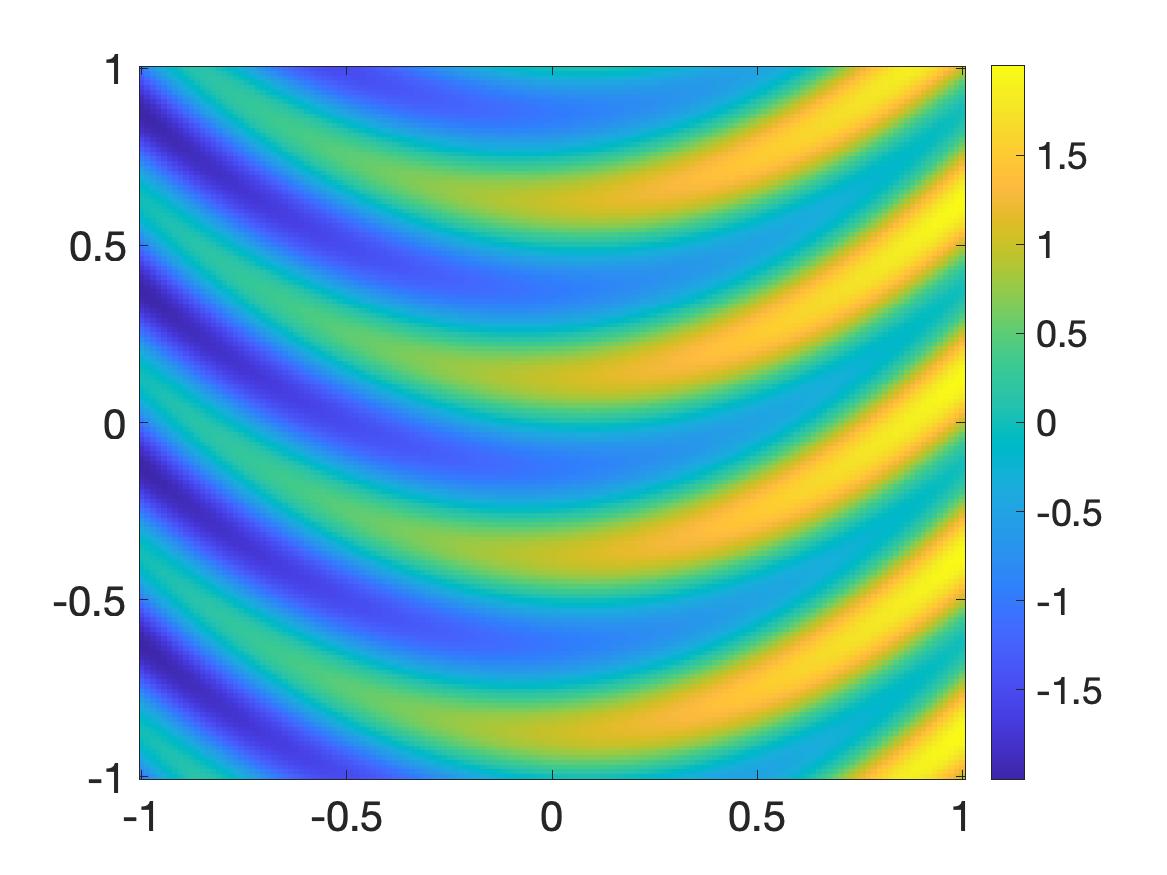}}
    \hfill
    \subfloat[The relative error $\frac{|u^* - u|}{\|u^*\|_{L^\infty(\Omega)}}$ when $\delta = 0\%$ ]{\includegraphics[width=.3\textwidth]{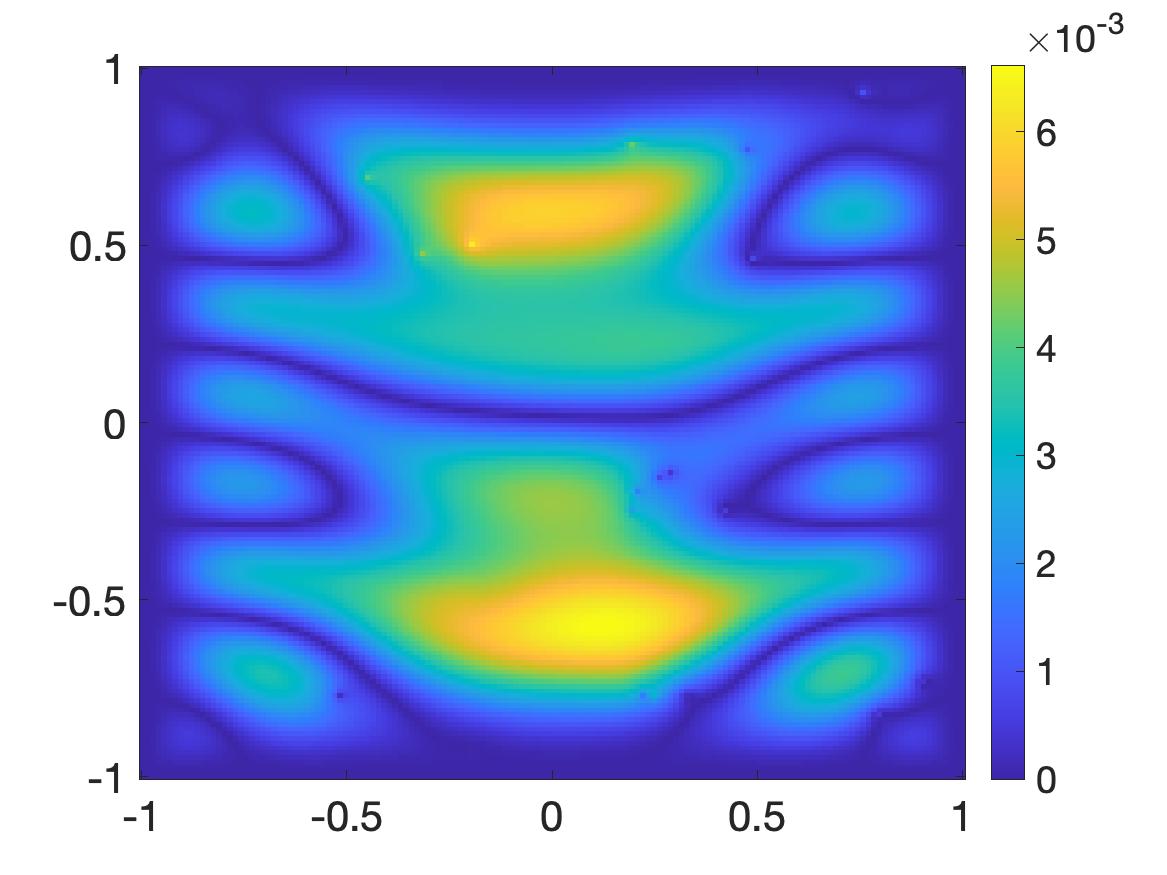}}
    
    \subfloat[The computed solution $u$  when $\delta = 10\%$ ]{\includegraphics[width=.3\textwidth]{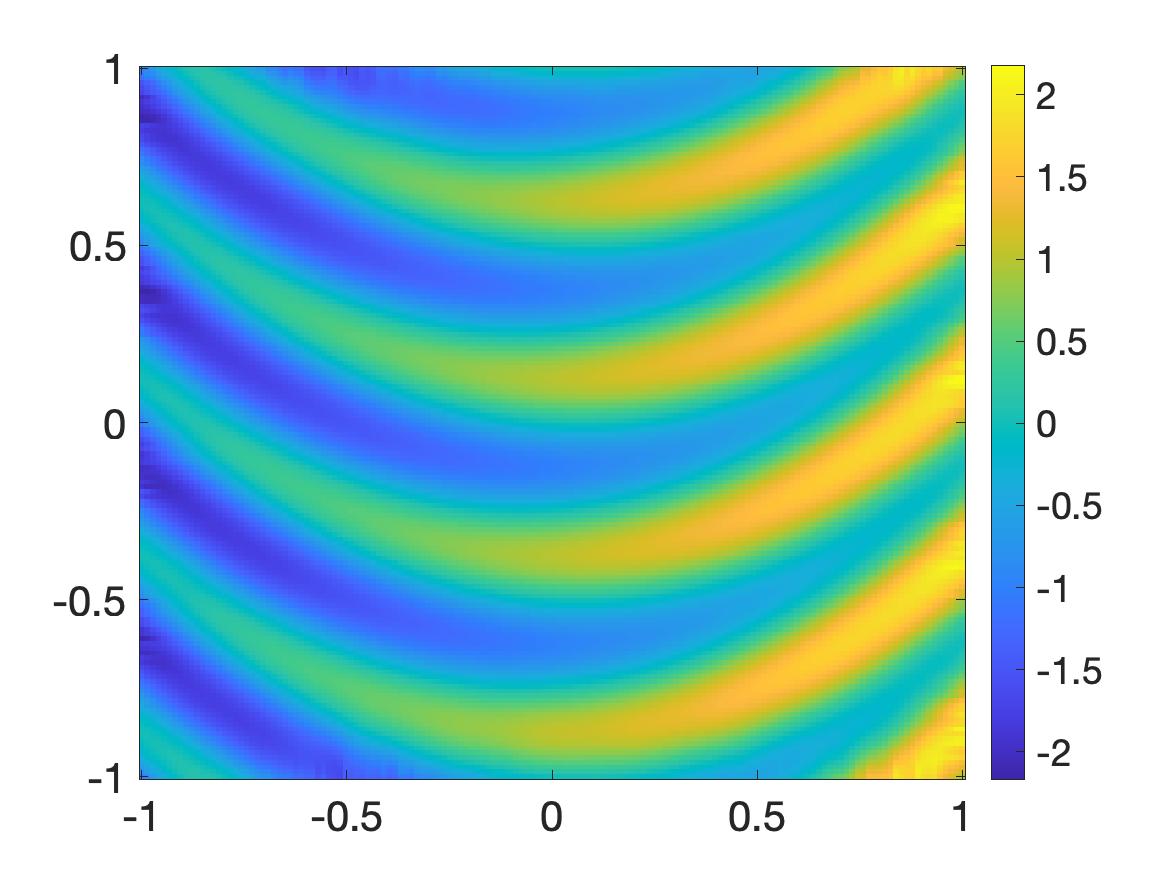}}
    \hfill
    \subfloat[\label{4e} The relative error $\frac{|u^* - u|}{\|u^*\|_{L^\infty(\Omega)}}$ when $\delta = 10\%$ ]{\includegraphics[width=.3\textwidth]{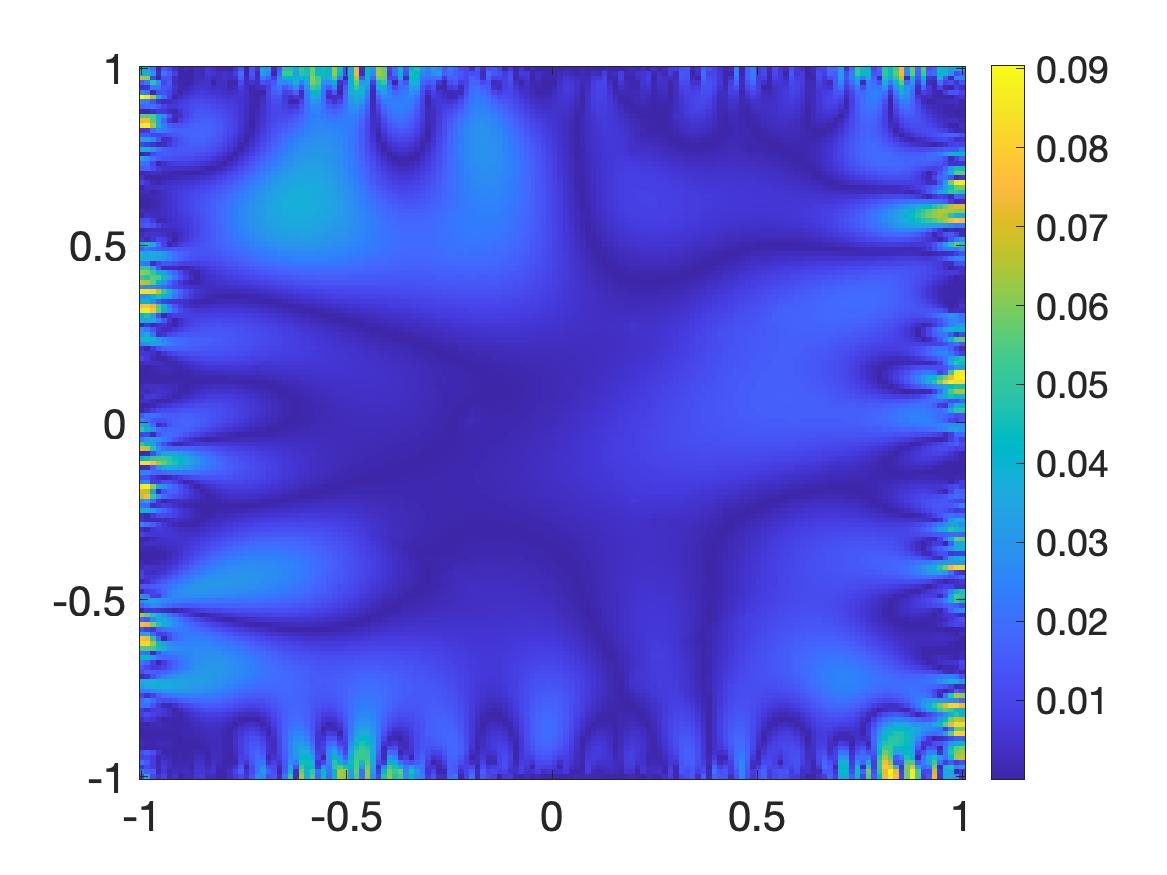}}
    \hfill
    \subfloat[\label{fig4f}The difference $\|u_{n+1} - u_n\|_{L^\infty(\Omega)}$ ]{\includegraphics[width=.3\textwidth]{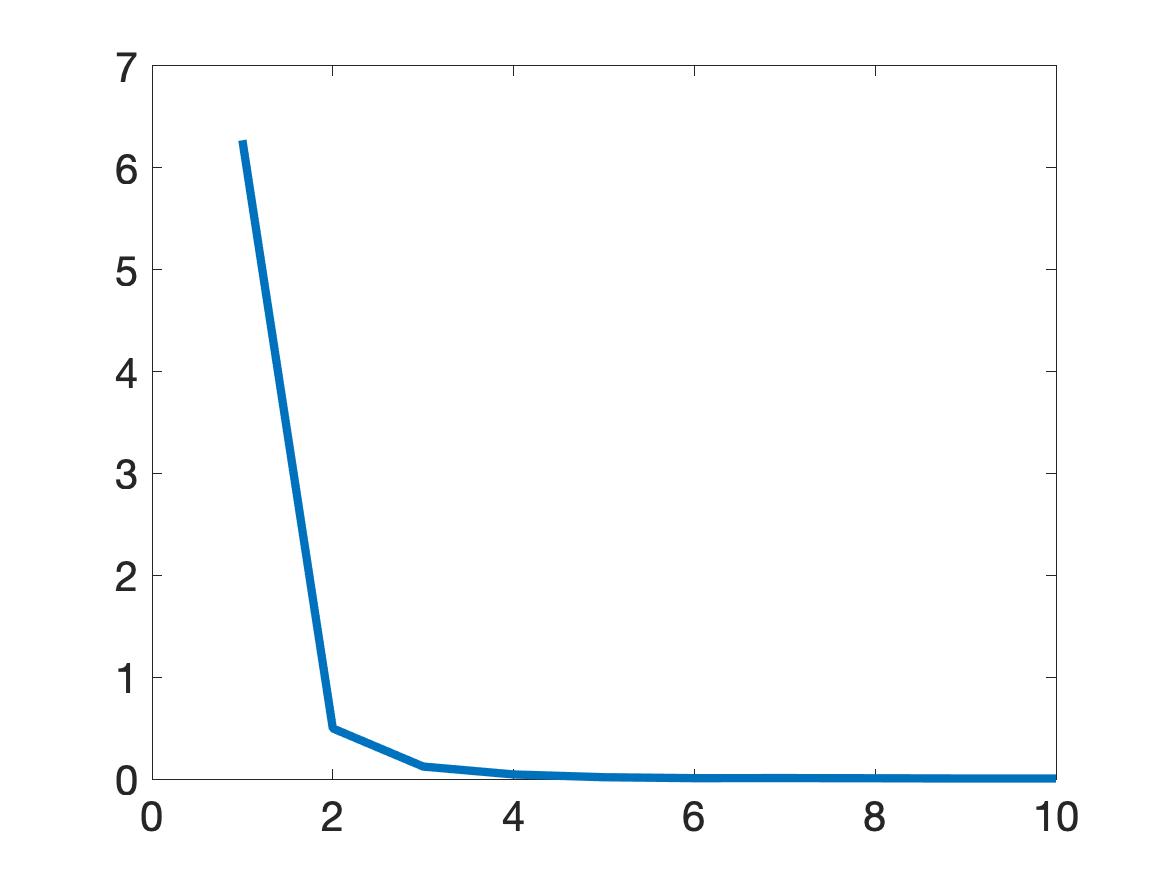}}
    \caption{Test 4. The graphs of the true and computed solution to \eqref{eqn4}, \eqref{f4} and \eqref{g4}  with noiseless and noisy boundary data. }
    \label{fig4}
\end{figure}

\noindent {\bf Test 4. }
We now test a more interesting problem when the nonlinearity $F(\x, s, p)$ grows at the quadratic rate in $s$ and is discontinuous with respect to $p$.
Let 
\[
    G(\x, s, p) = \left\{
        \begin{array}{ll}
             s^2 - e^{p_2}&  \mbox{if } e^{p_2} < 30,\\
             0& \mbox{otherwise}
        \end{array}
    \right.
\]
for all $\x \in \Omega,$ $s \in \R$, $p = (p_1, p_2) \in \R^2.$
We numerically solve the equation
\begin{equation}
    \Delta u(\x) + G(\x, u(\x), \nabla u(\x))
    -
    \Big[
       \big(\sin(4\pi x - 2\pi y^2) + y\big)^2
       -e^{-4\pi y \cos(4\pi x - 2\pi y^2) + 1}
    \Big] = 0
    \label{eqn4}
\end{equation}
for all $\x = (x, y) \in \Omega.$
The boundary data are given by
\begin{align}
    u(\x) &= \big(\sin(4\pi x - 2\pi y^2) + y\big)(1 + \delta\mbox{rand}_1), \label{f4}\\
    \partial_{\nu} u(\x)  &= 4\pi \big(
        \cos(4\pi x - 2\pi y^2),
        -y \cos(4\pi x - 2\pi y^2) + 1
    \big)\cdot \nu (1 + \delta \mbox{rand}_2) \label{g4}
\end{align}
for all $\x = (x, y) \in \partial \Omega$,
where $\delta > 0$ is the noise level and $\mbox{rand}_{i}$, $i = 1, 2$, is the function taking uniformly distributed random numbers in the rank $[-1, 1]$.
The true solution to \eqref{eqn4}, \eqref{f4} and \eqref{g4} when $\delta = 0$ is $u^*(\x) = \sin(4\pi x - 2\pi y^2) + y$ for all $\x = (x, y) \in \Omega.$
The graphs of the true and computed solution and the absolute error in the computation are displayed in Figure \ref{fig4}.

Even when \eqref{eqn4} involves a term that is not continuous with respect to $u_y$, Algorithm \ref{alg} delivers acceptable numerical solutions.
The relative errors in the computation are in Table \ref{table4}. On the other hand, one can observe from Figure \ref{fig4f} that our method converges fast.
The stopping criterion meets after only ten iterations. 

\begin{table}[h!]
    \centering
    \begin{tabular}{c c c }
    \hline
       Noise level  & $\|u^* - u^{\rm comp}\|_{L^\infty(\Omega)}$ & $\|u^* - u^{\rm comp}\|_{L^2(\Omega)}$ \\
      \hline
      $ \delta = 0\%$  
      &$0.0066$ & 
      $0.0055$ 
      \\
      $ \delta = 2\%$  
      &$0.0184$ & 
      $ 0.0104$ 
      \\
      $ \delta = 5\%$  
      &$ 0.0496$ & 
      $ 0.0179$ 
      \\
      $ \delta = 10\%$  
      &$0904$ & 
      $0.0325$ 
      \\
      \hline
    \end{tabular}
    \caption{Test 4. The relative errors in computation}
    \label{table4}
\end{table}

\begin{remark}
    We use a Macbook Pro 6-Core Intel Core i7 (2.6 GHz) to compute the numerical solutions above. The computational time for tests 1, 2, 3, and 4 are about 5 seconds, 7 seconds, 5 seconds, and 10 seconds, respectively. 
    The computational cost is not expensive.
\end{remark}

\begin{remark}
	In theory in and in the problem statement, we assume that we know an upper bound of $\|u^*\|_{C^1(\overline \Omega)}$, namely $M$, in order to use the cut-off technique in \eqref{cutoff}.
	However, in the numerical study, this step was not implemented.
	That means we relax the request that we need to know $M$. 
	The method might be stronger than what we can rigorously prove.
\end{remark}

\begin{remark}
The numerical method in Algorithm \ref{alg} is stronger than what we can prove in Theorem \ref{thm1} and Theorem \ref{thm2}. It can deliver reliable solutions even when the nonlinearity is not smooth (see test 3) and is not continuous (see test 4). 
\end{remark}

%

%

%
%
%
%
%
%
%
%
%
%
\section{Concluding remarks}\label{sec6}
We have solved the problem of computing solutions to quasi-linear PDEs. 
Although this problem is nonlinear,  we do not require a good initial guess of the true solution.
We first define an operator $\Phi$ such that the true solution to the given quasilinear PDE is the fixed point of $\Phi$.
We construct a recursive sequence $\{u_n\}_{n \geq 0}$ whose initial term $u_0$ can be taken arbitrary and the $n^{\rm th}$ term $u_n = \Phi(u_{n - 1})$.
We next apply a Carleman estimate to prove the convergence of this sequence.
Moreover, we have proved that the stability of our method with respect to noise is of the Lipschitz type.
Some interesting numerical examples are presented.

\section*{Acknowledgement}   
This work is dedicated to Professor Duong Minh Duc.
It was partially supported by National Science Foundation grant DMS-2208159 and by funds
provided by the Faculty Research Grant program at UNC Charlotte Fund No. 11127


\end{document}